\documentstyle[12pt,amstex,epic,psfig,epsfig]{article}

\catcode`\@=11
\def\numberbysection{\@addtoreset{equation}{section}
        \def\theequation{\thesection.\arabic{equation}}}
\numberbysection
\newtheorem{theo}{Theorem}
\newtheorem{lemma}{Lemma}
\newtheorem{prop}{Proposition}

\newtheorem{defi}{Definition}
\def\N{{\cal N}}
\def\=p{ = {\hspace{-2.6 ex}\raisebox{-1.1 ex}{\scriptsize{$(p)$}}}}
\def\S{{\cal S}}
\def\I{{\cal I}}

\begin{document}

\baselineskip 18pt

\centerline{\LARGE  From $sl(2)$ Kirby weight sytems }

\vskip .5cm

\centerline{\LARGE to the asymptotic 3-manifold invariant }

\vskip 1cm
\centerline{\large \bf Laurent Freidel
\footnote{On leave of ENSLAPP,ENS-Lyon, FRANCE.  This work is supported
by CNRS and  NATO grant.\\
email:freidel@enslapp-ens-lyon.fr}}

\vskip 1cm
\centerline{\bf
Center for Gravitational Physics and Geometry}

\centerline{\bf
Penn State University }

\centerline{\bf U.S.A}

\centerline{\bf
e-mail:freidel@phys.psu.edu }

\vskip 12pt
\centerline{}

\vskip 1cm

\begin{abstract}
We give a construction of  Kirby weight systems associated to $sl(2)$
and valued into the finite field  ${{\Bbb Z} / p{\Bbb Z}}$.
We show that it is possible to apply this sequence of 
 weight systems on the universal invariant of framed link.
 We also show that 
the corresponding sequence admits a Fermat limit, which defines an 
asymptotic rational homology 3-sphere quantum invariant.  Moreover, 
this asymptotic invariant coincides with the Ohtsuki invariant.
\end{abstract}
\newpage

\section{Introduction}
The universal invariant for oriented framed link 
 is now understood to be the most fundamental object 
in the study of quantum invariant of knots and links.
One of the main features of this invariant is that it contains, at 
least, the same information as all the Witten-Reshetikhin-Turaev 
quantum invariants associated with semi-simple quantum groups 
\cite{RT}.  There exist by now three constructive definition of this 
invariant~: The first one was given by Kontsevich using the 
Khniznik-Zamolodchikov equation\cite{Kon}, 
the second was given by Piunikine \cite{Piu} and Cartier \cite{Car} using a 
combinatorial definition inspired by the quasi-Hopf algebras 
construction of Drinfeld \cite{Dri} and the third was given by the 
perturbative expansion of Chern-Simons theory (\cite{AF1,BN}).  The 
first two have been proved to be equal \cite{AF1,LeM1} and are called 
canonical, but to my knowledge it is still a open question to show that 
the third one is also canonical.

Using the surgery presentation of a 3-dimensional manifold a framed link
 invariant can be promoted to a 3-d manifold invariant if it proves
 to be invariant under Kirby moves.
Using the representation theory of quantum groups at roots of unity, 
Reshetikhin and Turaev constructed a quantum invariant of 3-manifolds 
as a weighted sum of framed colored Jones polynomials (in the case of 
$SU(2)$) over finite dimensional irreducible highest weight 
representations \cite{RT2}.  This formula was conjectured by Witten to 
be the partition function of Chern-Simons theory \cite{Wit}.  This 
construction uses the fact that representation theory of quantum groups 
at roots of unity is truncated, so the weighted sum involves only a 
finite number of terms.  The restriction to the case of roots of unity 
appears as a peculiar regularisation of Chern-Simons theory.  But it 
is also known from general grounds that behind every regularisation of a 
sensible quantum field theory there exists a renormalized theory which is 
regularisation independent.  In this context this means that it should 
be possible to define an asymptotic partition function for 
Chern-Simons as a formal power series in $\hbar$.  The quest for such 
a theory has been concentrated on the value of the Chern-Simons theory for 
homology sphere expanded around the trivial connection, which is  an 
isolated critical point of the Chern-Simons action.  The oldest 
technique at disposal is the usual perturbative definition of quantum 
field theory \cite{AxS} but it was only very recently that ther appeared the 
possibility of defining Chern-Simons theory as a perturbative expansion 
around the trivial connection \cite{Bot}.  Roszansky developed 
techniques in order to study the asymptotic behaviour of W-R-T 
invariants in the limit where $k$ (the order of the root of unity) goes 
to infinity.  In this context Ohtsuki \cite{Oht} showed that the 
$su(2)$ W-R-T for integral and rational homology spheres admits
a "Fermat limit" (we use the terminology of \cite{Lin}) and thus gives rise to $su(2)$ homology sphere invariants for 
generic $q$.  The breakthrough in this quest was achieved in three 
seminal papers~: \cite{LeM2,LeM3,Le2}, where 
 T.Q.~Le, H.~Murakami, J.~Murakami,  T.~Ohtsuki, first showed that it 
was possible to construct a 3-manifold invariant from the universal 
framed link invariant if one could find a weight system satisfying the 
so called Kirby relations (we call such a weight system a Kirby weight 
system).  Then they constructed a universal Kirby weight system which 
maps the space of chord diagrams to the space of trivalent graphs 
satisfying Jacobi (or so called I-H-X) relations.  Moreover, T.Q.~Le 
showed \cite{Le2} that the resulting object obtained by applying this Kirby 
weight system is the universal invariant of rational homology sphere, 
i-e the graded isomorphisms between Jacobi trivalent graphs and finite 
type homology spheres invariants.  This result is the 3-manifold 
analogue of the celebrated theorem of Kontsevich for knots \cite{Kon}.  
One can think that this is the end of the story.  But a lot of 
questions were remaining like the correspondence between this 
universal invariants and the numerical ones obtained from the W-R-T or 
Ohtsuki constructions.  Of course, Ohtsuki \cite{Oht3} proved that by 
applying the $su(2)$ trivalent graphs weight system on the universal 
Homology sphere invariant one obtains the Ohtsuki invariant, see also 
theorem 8.5 of \cite{Oht2}.

Our aim in this paper is to give a completely different 
construction of 
a Kirby weight system in the case of $su(2)$ and to evaluate the 
corresponding 3-manifold invariant.  Our construction is very 
reminiscent of the Reshetikhin Turaev construction in the sense that 
the the $su(2)$ Kirby weight system on chord diagrams we propose is 
obtained as a weighted sum of weight systems over irreducible 
representations.  In order to give a meaning to this sum we were forced 
to work in non zero characteristic fields.  this restriction appears 
analogous to roots of unity regularisation at the level of weight 
systems.  Using results of Le \cite{Le3} on integrality properties of 
the Universal framed link invariants we apply this weight system on 
this invariant.  We show that the corresponding object admits a 
"Fermat limit"  which is the 
Ohtsuki invariant given in terms of Gaussian integrals.  This fact 
appears as a specific realization in the case of $su(2)$ of the 
conjectural proposal given in \cite{BN1}.  Moreover this limit is the 
same for $su(2)$ and $so(3)$ invariants.  The first part of this paper 
was presented in Knots 96 Tokyo conference.  The study of the general 
Lie algebra case is under construction in collaboration with D.  
Altschuler.

In section 2 we present the construction of skein modules in non zero 
characteristic fields.
The material of this part is largely borrowed from \cite{turaevb}
 (see also \cite{Kaufflin}),
 but I couldn't find in the literature any references concerning
 the finite field characteristic case (even if it is very similar to the root
 of unity case).  I also present in this section 
 some recoupling theory 
 theorem that I uses in the following.  The third section is devoted to 
 the construction of weight system from skein modules, the link between 
 zero and non zero characteristic field weight system and the link with 
 the usual definition of $su(2)$ weight system.  The section 4 presents 
 the construction of a $so(3)$ and $su(2)$ weight system, the proof 
 that they satisfy the Kirby relations and the study of some of the important 
 properties that these weight systems satisfy.  In the section 5 we take 
 advantage of these properties in order to apply the Kirby weight 
 systems the Universal Framed link invariant.  We then show that the 
 resulting object admits a Fermat limit expressed in terms of Gaussian 
 integrals.  This fermat limit is the same for the $su(2)$ and $ so(3)$ 
 case and moreover they coincide with Ohtsuki construction.

\section{Skein Module}

{\bf Preliminaries} \\
In this article $p$ will denote a prime odd integer.
$ {{\Bbb{Z}} / p {\Bbb{Z}}} $ is identified with the set 
$\{0,\cdots,p-1\}$ and the set of invertible element 
of $ {{\Bbb{Z}} / p {\Bbb{Z}}} $ with $\{1,\cdots,p-1\}$.
If $k \in {\Bbb{Z}}$ we denote by $\psi_p(k)$ its image in ${{\Bbb{Z}} 
/ p {\Bbb{Z}}}$; moreover, if $k$ is not divisible by $p$ then 
$\psi_p(k)$ is invertible in $ {{\Bbb{Z}} / p {\Bbb{Z}}} $ and we 
denotes its inverse by $\psi_p({1\over k})$.  Lettting ${\Bbb{Q}}_p=\{ 
{a\over b} \in {\Bbb{Q}},\,gcd(a,b)=1,\, gcd(p,b)=1\} $, we define the 
homomorphism $\psi_p: {\Bbb{Q}}_p \rightarrow {{\Bbb{Z}} / p 
{\Bbb{Z}}}$ by $\psi_p({a\over b} )={\psi_p(a) \psi_p({1\over b})}$.  
If $x,y \in {\Bbb{Q}}_p $, $x \=p y$ means $\psi_p(x) =\psi_p(y)$.  \\

\begin{figure}[htbp]
\centerline{
\begin{picture}(0,0)%
\epsfig{file=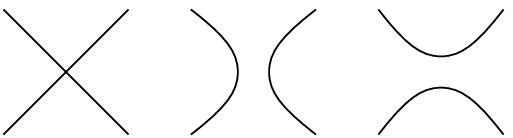}%
\end{picture}%
\setlength{\unitlength}{0.00083300in}%
\begingroup\makeatletter\ifx\SetFigFont\undefined
\def\x#1#2#3#4#5#6#7\relax{\def\x{#1#2#3#4#5#6}}%
\expandafter\x\fmtname xxxxxx\relax \def\y{splain}%
\ifx\x\y   
\gdef\SetFigFont#1#2#3{%
  \ifnum #1<17\tiny\else \ifnum #1<20\small\else
  \ifnum #1<24\normalsize\else \ifnum #1<29\large\else
  \ifnum #1<34\Large\else \ifnum #1<41\LARGE\else
     \huge\fi\fi\fi\fi\fi\fi
  \csname #3\endcsname}%
\else
\gdef\SetFigFont#1#2#3{\begingroup
  \count@#1\relax \ifnum 25<\count@\count@25\fi
  \def\x{\endgroup\@setsize\SetFigFont{#2pt}}%
  \expandafter\x
    \csname \romannumeral\the\count@ pt\expandafter\endcsname
    \csname @\romannumeral\the\count@ pt\endcsname
  \csname #3\endcsname}%
\fi
\fi\endgroup
\begin{picture}(2427,624)(2089,-1573)
\put(2708,-1296){\makebox(0,0)[lb]{\smash{\SetFigFont{12}{14.4}{rm}
$+$
}}}
\put(3615,-1303){\makebox(0,0)[lb]{\smash{\SetFigFont{12}{14.4}{rm}
$+$
}}}
\put(4516,-1288){\makebox(0,0)[lb]{\smash{\SetFigFont{12}{14.4}{rm}
$=$ $0$
}}}
\end{picture}
 }
\caption[]{}
\label{skein}
\end{figure}

\begin{figure}[htbp]
\centerline{
\begin{picture}(0,0)%
\epsfig{file=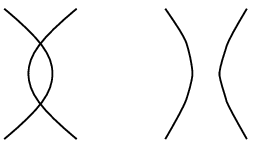}%
\end{picture}%
\setlength{\unitlength}{0.00083300in}%
\begingroup\makeatletter\ifx\SetFigFont\undefined
\def\x#1#2#3#4#5#6#7\relax{\def\x{#1#2#3#4#5#6}}%
\expandafter\x\fmtname xxxxxx\relax \def\y{splain}%
\ifx\x\y   
\gdef\SetFigFont#1#2#3{%
  \ifnum #1<17\tiny\else \ifnum #1<20\small\else
  \ifnum #1<24\normalsize\else \ifnum #1<29\large\else
  \ifnum #1<34\Large\else \ifnum #1<41\LARGE\else
     \huge\fi\fi\fi\fi\fi\fi
  \csname #3\endcsname}%
\else
\gdef\SetFigFont#1#2#3{\begingroup
  \count@#1\relax \ifnum 25<\count@\count@25\fi
  \def\x{\endgroup\@setsize\SetFigFont{#2pt}}%
  \expandafter\x
    \csname \romannumeral\the\count@ pt\expandafter\endcsname
    \csname @\romannumeral\the\count@ pt\endcsname
  \csname #3\endcsname}%
\fi
\fi\endgroup
\begin{picture}(1180,641)(3147,-1545)
\put(3601,-1261){\makebox(0,0)[lb]{\smash{\SetFigFont{12}{14.4}{rm}
$=$
}}}
\end{picture}
 }
\caption[]{}
\label{id}
\end{figure}

\begin{figure}[htbp]
\centerline{
\begin{picture}(0,0)%
\epsfig{file=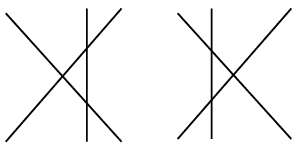}%
\end{picture}%
\setlength{\unitlength}{0.00083300in}%
\begingroup\makeatletter\ifx\SetFigFont\undefined
\def\x#1#2#3#4#5#6#7\relax{\def\x{#1#2#3#4#5#6}}%
\expandafter\x\fmtname xxxxxx\relax \def\y{splain}%
\ifx\x\y   
\gdef\SetFigFont#1#2#3{%
  \ifnum #1<17\tiny\else \ifnum #1<20\small\else
  \ifnum #1<24\normalsize\else \ifnum #1<29\large\else
  \ifnum #1<34\Large\else \ifnum #1<41\LARGE\else
     \huge\fi\fi\fi\fi\fi\fi
  \csname #3\endcsname}%
\else
\gdef\SetFigFont#1#2#3{\begingroup
  \count@#1\relax \ifnum 25<\count@\count@25\fi
  \def\x{\endgroup\@setsize\SetFigFont{#2pt}}%
  \expandafter\x
    \csname \romannumeral\the\count@ pt\expandafter\endcsname
    \csname @\romannumeral\the\count@ pt\endcsname
  \csname #3\endcsname}%
\fi
\fi\endgroup
\begin{picture}(1394,663)(2717,-2814)
\put(3271,-2536){\makebox(0,0)[lb]{\smash{\SetFigFont{12}{14.4}{rm}
$=$
}}}
\end{picture}
 }
\caption[]{}
\label{YB}
\end{figure}

\begin{figure}[htbp]
\centerline{
\begin{picture}(0,0)%
\epsfig{file=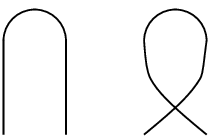}%
\end{picture}%
\setlength{\unitlength}{0.00083300in}%
\begingroup\makeatletter\ifx\SetFigFont\undefined
\def\x#1#2#3#4#5#6#7\relax{\def\x{#1#2#3#4#5#6}}%
\expandafter\x\fmtname xxxxxx\relax \def\y{splain}%
\ifx\x\y   
\gdef\SetFigFont#1#2#3{%
  \ifnum #1<17\tiny\else \ifnum #1<20\small\else
  \ifnum #1<24\normalsize\else \ifnum #1<29\large\else
  \ifnum #1<34\Large\else \ifnum #1<41\LARGE\else
     \huge\fi\fi\fi\fi\fi\fi
  \csname #3\endcsname}%
\else
\gdef\SetFigFont#1#2#3{\begingroup
  \count@#1\relax \ifnum 25<\count@\count@25\fi
  \def\x{\endgroup\@setsize\SetFigFont{#2pt}}%
  \expandafter\x
    \csname \romannumeral\the\count@ pt\expandafter\endcsname
    \csname @\romannumeral\the\count@ pt\endcsname
  \csname #3\endcsname}%
\fi
\fi\endgroup
\begin{picture}(996,621)(2089,-1273)
\put(2476,-1041){\makebox(0,0)[lb]{\smash{\SetFigFont{12}{14.4}{rm}
$=$
}}}
\end{picture}
 }
\caption[]{}
\label{writ}
\end{figure}

A $(k,l)$ Tangle $T$ is an immersion of a one-dimensional compact 
submanifold, considered modulo ambient isotopy of ${\Bbb{R}} \times 
[0,1]$, into ${\Bbb{R}} \times [0,1] $, such that $\partial T = T 
\bigcap \left[ {\Bbb{R}} \times \{0\} \cup \{1\}\right]$, and $ Card( T \cap 
[{\Bbb{R}} \times \{0\}] ) = k $, $ Card( T \cap [{\Bbb{R}} \times \{1\}] 
) = l $.  Denote ${\cal T} (k,l) $ the ${\Bbb{Z}}$-linear span of all 
$(k,l)$ tangles.
 \begin{defi}
We denote by ${\cal S}_p(k,l)$ the ${{\Bbb{Z}} / p {\Bbb{Z}}}$-module generated
 by $(k,l)$ tangles
with the following relations~:
\begin{eqnarray}
&T \cup O  = -2 \, T  \label{circ} \mbox{ where $O$ denote the circle 
and $T$ is an arbitrary tangle},\\
&\label{skei} \mbox{the skein relation shown in figure \ref{skein}}.
\end{eqnarray}
\end{defi}
Lets us call simple the tangle diagrams without any crossing and 
without any circle; using the skein relation, we have the following 
lemma~:
\begin{lemma}
\label{isotopy}
${\cal S}_p{(k,l)}$ is a free $K_p$ module with basis given by simple $(k,l)$ 
diagrams. In ${\cal S}_p(k,l)$ the relations shown in
figures \ref{id}, \ref{YB}, \ref{writ} are satisfied~:
\end{lemma}

Using the skein property any tangle can be decomposed uniquely into a 
sum of diagrams without any crossings.
Using \ref{circ}, this sum can be expanded into a sum of simple 
diagrams.  Thus we get the first property of the lemma.  The second 
part of the lemma is obtained by direct computations. In particular, the 
basis of ${\cal S}_p(0,0)$ is given by the empty tangle, so this gives 
an identification of ${\cal S}_p{(0,0)}$ with ${{\Bbb{Z}} / p {\Bbb{Z}}}$.  
If $x\in {\cal S}_p(k,n)$ and $ y\in {\cal S}_p(l,k)$ we can define 
the composition $x\cdot y \in {\cal S}_p(l,n)$ by stacking elements 
(see fig \ref{stac}) and the tensor product$x\otimes y$ by juxtaposing elements 
(see fig \ref{tens}).
\begin{figure}[htbp]
\centerline{
\setlength{\unitlength}{0.0125in}
\begingroup\makeatletter\ifx\SetFigFont\undefined
\def\x#1#2#3#4#5#6#7\relax{\def\x{#1#2#3#4#5#6}}%
\expandafter\x\fmtname xxxxxx\relax \def\y{splain}%
\ifx\x\y   
\gdef\SetFigFont#1#2#3{%
  \ifnum #1<17\tiny\else \ifnum #1<20\small\else
  \ifnum #1<24\normalsize\else \ifnum #1<29\large\else
  \ifnum #1<34\Large\else \ifnum #1<41\LARGE\else
     \huge\fi\fi\fi\fi\fi\fi
  \csname #3\endcsname}%
\else
\gdef\SetFigFont#1#2#3{\begingroup
  \count@#1\relax \ifnum 25<\count@\count@25\fi
  \def\x{\endgroup\@setsize\SetFigFont{#2pt}}%
  \expandafter\x
    \csname \romannumeral\the\count@ pt\expandafter\endcsname
    \csname @\romannumeral\the\count@ pt\endcsname
  \csname #3\endcsname}%
\fi
\fi\endgroup
\begin{picture}(235,125)(0,-10)
\drawline(90,5)(90,45)
\drawline(130,5)(130,45)
\drawline(85,65)(85,45)(135,45)
	(135,65)(85,65)(85,65)
\put(107,51){\makebox(0,0)[lb]{\smash{{{\SetFigFont{12}{14.4}{rm}g}}}}}
\drawline(0,65)(0,45)(50,45)
	(50,65)(0,65)(0,65)
\drawline(5,5)(5,45)
\drawline(45,5)(45,45)
\drawline(45,105)(45,65)
\drawline(5,105)(5,65)
\put(23,52){\makebox(0,0)[lb]{\smash{{{\SetFigFont{12}{14.4}{rm}f}}}}}
\drawline(185,88)(185,68)(235,68)
	(235,88)(185,88)(185,88)
\drawline(230,45)(230,68)
\drawline(230,110)(230,88)
\drawline(190,110)(190,88)
\drawline(190,45)(190,68)
\drawline(185,44)(185,24)(235,24)
	(235,44)(185,44)(185,44)
\drawline(190,0)(190,24)
\drawline(230,0)(230,24)
\put(208,75){\makebox(0,0)[lb]{\smash{{{\SetFigFont{12}{14.4}{rm}f}}}}}
\put(207,30){\makebox(0,0)[lb]{\smash{{{\SetFigFont{12}{14.4}{rm}g}}}}}
\put(67,55){\circle{6}}
\drawline(130,105)(130,65)
\drawline(90,105)(90,65)
\put(152,51){\makebox(0,0)[lb]{\smash{{{\SetFigFont{12}{14.4}{rm}=}}}}}
\end{picture}
 }
\caption[]{}
\label{stac}
\end{figure}
\begin{figure}[htbp]
\centerline{
\setlength{\unitlength}{0.0125in}
\begingroup\makeatletter\ifx\SetFigFont\undefined
\def\x#1#2#3#4#5#6#7\relax{\def\x{#1#2#3#4#5#6}}%
\expandafter\x\fmtname xxxxxx\relax \def\y{splain}%
\ifx\x\y   
\gdef\SetFigFont#1#2#3{%
  \ifnum #1<17\tiny\else \ifnum #1<20\small\else
  \ifnum #1<24\normalsize\else \ifnum #1<29\large\else
  \ifnum #1<34\Large\else \ifnum #1<41\LARGE\else
     \huge\fi\fi\fi\fi\fi\fi
  \csname #3\endcsname}%
\else
\gdef\SetFigFont#1#2#3{\begingroup
  \count@#1\relax \ifnum 25<\count@\count@25\fi
  \def\x{\endgroup\@setsize\SetFigFont{#2pt}}%
  \expandafter\x
    \csname \romannumeral\the\count@ pt\expandafter\endcsname
    \csname @\romannumeral\the\count@ pt\endcsname
  \csname #3\endcsname}%
\fi
\fi\endgroup
\begin{picture}(265,115)(0,-10)
\drawline(0,60)(0,40)(50,40)
	(50,60)(0,60)(0,60)
\drawline(5,0)(5,40)
\drawline(45,0)(45,40)
\drawline(45,100)(45,60)
\drawline(5,100)(5,60)
\put(22,45){\makebox(0,0)[lb]{$f$}}
\drawline(165,60)(165,40)(215,40)
	(215,60)(165,60)(165,60)
\drawline(170,0)(170,40)
\drawline(210,0)(210,40)
\drawline(210,100)(210,60)
\drawline(170,100)(170,60)
\put(188,47){\makebox(0,0)[lb]{\smash{{{\SetFigFont{12}{14.4}{it}f}}}}}
\drawline(90,0)(90,40)
\drawline(130,0)(130,40)
\drawline(85,60)(85,40)(135,40)
	(135,60)(85,60)(85,60)
\put(107,46){\makebox(0,0)[lb]{\smash{{{\SetFigFont{12}{14.4}{rm}g}}}}}
\drawline(130,100)(130,60)
\drawline(90,100)(90,60)
\drawline(220,0)(220,40)
\drawline(260,0)(260,40)
\drawline(215,60)(215,40)(265,40)
	(265,60)(215,60)(215,60)
\put(237,46){\makebox(0,0)[lb]{\smash{{{\SetFigFont{12}{14.4}{rm}g}}}}}
\drawline(260,100)(260,60)
\drawline(220,100)(220,60)
\put(63,48){$\otimes$}
\put(145,45){\makebox(0,0)[lb]{\smash{{{\SetFigFont{12}{14.4}{rm}$=$}}}}}
\end{picture}
 }
\caption[]{}
\label{tens}
\end{figure}
We also define a trace $tr f$ to be the closure of f as in figure \ref{tr}.
In ${\cal S}_p{(k,l)}$ we consider the submodule of null elements $\N_p{(k,l)}$ 
defined as follows :
\begin{equation}
\N_p{(k,l)} =\left\{ n \in {\cal S}_p{(k,l)} | \forall f \in 
{\cal S}_p{(l,k)},\, 
tr( fn)\=p 0 \right\}
\end{equation}

\begin{figure}[htbp]
\centerline{
\setlength{\unitlength}{0.0125in}
\begingroup\makeatletter\ifx\SetFigFont\undefined
\def\x#1#2#3#4#5#6#7\relax{\def\x{#1#2#3#4#5#6}}%
\expandafter\x\fmtname xxxxxx\relax \def\y{splain}%
\ifx\x\y   
\gdef\SetFigFont#1#2#3{%
  \ifnum #1<17\tiny\else \ifnum #1<20\small\else
  \ifnum #1<24\normalsize\else \ifnum #1<29\large\else
  \ifnum #1<34\Large\else \ifnum #1<41\LARGE\else
     \huge\fi\fi\fi\fi\fi\fi
  \csname #3\endcsname}%
\else
\gdef\SetFigFont#1#2#3{\begingroup
  \count@#1\relax \ifnum 25<\count@\count@25\fi
  \def\x{\endgroup\@setsize\SetFigFont{#2pt}}%
  \expandafter\x
    \csname \romannumeral\the\count@ pt\expandafter\endcsname
    \csname @\romannumeral\the\count@ pt\endcsname
  \csname #3\endcsname}%
\fi
\fi\endgroup
\begin{picture}(209,115)(0,-10)
\drawline(25,60)(25,40)(75,40)
	(75,60)(25,60)(25,60)
\drawline(30,0)(30,40)
\drawline(70,0)(70,40)
\drawline(70,100)(70,60)
\drawline(30,100)(30,60)
\put(48,47){\makebox(0,0)[lb]{\smash{{{\SetFigFont{12}{14.4}{rm}f}}}}}
\drawline(110,60)(110,40)(160,40)
	(160,60)(110,60)(110,60)
\drawline(155,40)	(157.582,35.121)
	(159.953,31.984)
	(162.347,30.355)
	(165.000,30.000)
\drawline(165,30)	(168.224,30.696)
	(171.135,32.207)
	(175.870,37.073)
	(177.623,40.129)
	(178.920,43.402)
	(179.724,46.742)
	(180.000,50.000)
\drawline(180,50)	(179.724,53.258)
	(178.920,56.598)
	(177.623,59.871)
	(175.870,62.927)
	(171.135,67.793)
	(168.224,69.304)
	(165.000,70.000)
\drawline(165,70)	(162.347,69.645)
	(159.953,68.016)
	(157.582,64.879)
	(155.000,60.000)
\drawline(115,40)	(117.888,35.325)
	(120.655,31.042)
	(123.317,27.137)
	(125.889,23.598)
	(128.384,20.412)
	(130.818,17.566)
	(135.560,12.844)
	(140.232,9.329)
	(144.951,6.918)
	(149.835,5.509)
	(152.375,5.149)
	(155.000,5.000)
\drawline(155,5)	(159.701,5.145)
	(164.433,5.729)
	(169.152,6.734)
	(173.814,8.139)
	(178.374,9.925)
	(182.787,12.071)
	(187.008,14.558)
	(190.993,17.366)
	(194.697,20.475)
	(198.075,23.865)
	(201.082,27.517)
	(203.675,31.410)
	(205.808,35.525)
	(207.436,39.841)
	(208.515,44.339)
	(209.000,49.000)
\drawline(209,49)	(209.013,51.517)
	(208.869,54.015)
	(208.130,58.939)
	(206.822,63.741)
	(204.982,68.386)
	(202.647,72.843)
	(199.855,77.078)
	(196.642,81.058)
	(193.047,84.751)
	(189.107,88.123)
	(184.859,91.142)
	(180.340,93.775)
	(175.588,95.990)
	(170.641,97.752)
	(168.106,98.454)
	(165.536,99.030)
	(162.936,99.477)
	(160.310,99.790)
	(157.663,99.966)
	(155.000,100.000)
\drawline(155,100)	(152.021,99.840)
	(149.175,99.436)
	(146.446,98.773)
	(143.820,97.838)
	(141.282,96.616)
	(138.818,95.091)
	(134.052,91.076)
	(131.721,88.557)
	(129.404,85.677)
	(127.088,82.421)
	(124.758,78.776)
	(122.398,74.726)
	(119.995,70.256)
	(117.534,65.352)
	(116.277,62.733)
	(115.000,60.000)
\put(0,45){\makebox(0,0)[lb]{\smash{{{\SetFigFont{12}{14.4}{rm}tr}}}}}
\put(90,45){\makebox(0,0)[lb]{\smash{{{\SetFigFont{12}{14.4}{rm}=}}}}}
\put(133,47){\makebox(0,0)[lb]{\smash{{{\SetFigFont{12}{14.4}{rm}f}}}}}
\end{picture}
 }
\caption[]{}
\label{tr}
\end{figure}
\begin{figure}[htbp]
\centerline{
\begin{picture}(0,0)%
\epsfig{file=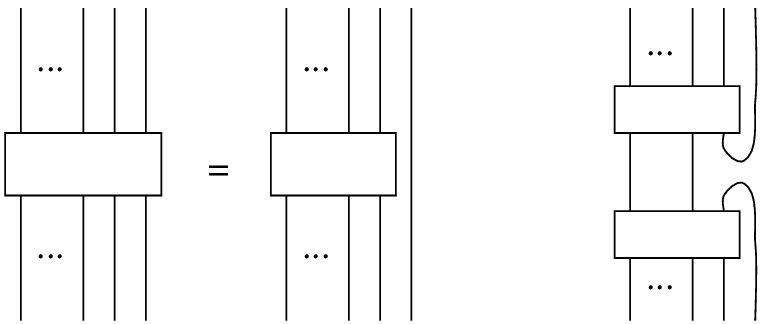}%
\end{picture}%
\setlength{\unitlength}{0.00083300in}%
\begingroup\makeatletter\ifx\SetFigFont\undefined
\def\x#1#2#3#4#5#6#7\relax{\def\x{#1#2#3#4#5#6}}%
\expandafter\x\fmtname xxxxxx\relax \def\y{splain}%
\ifx\x\y   
\gdef\SetFigFont#1#2#3{%
  \ifnum #1<17\tiny\else \ifnum #1<20\small\else
  \ifnum #1<24\normalsize\else \ifnum #1<29\large\else
  \ifnum #1<34\Large\else \ifnum #1<41\LARGE\else
     \huge\fi\fi\fi\fi\fi\fi
  \csname #3\endcsname}%
\else
\gdef\SetFigFont#1#2#3{\begingroup
  \count@#1\relax \ifnum 25<\count@\count@25\fi
  \def\x{\endgroup\@setsize\SetFigFont{#2pt}}%
  \expandafter\x
    \csname \romannumeral\the\count@ pt\expandafter\endcsname
    \csname @\romannumeral\the\count@ pt\endcsname
  \csname #3\endcsname}%
\fi
\fi\endgroup
\begin{picture}(3626,1524)(1414,-3673)
\put(2971,-2956){\makebox(0,0)[lb]{\smash{\SetFigFont{12}{14.4}{rm}$f_n$}}}
\put(1771,-2956){\makebox(0,0)[lb]{\smash{\SetFigFont{12}{14.4}{rm}$f_{n+1}$}}}
\put(3601,-2986){\makebox(0,0)[lb]{\smash{\SetFigFont{12}{14.4}{rm}$+$}}}
\put(4606,-2701){\makebox(0,0)[lb]{\smash{\SetFigFont{12}{14.4}{rm}$f_n$}}}
\put(4606,-3301){\makebox(0,0)[lb]{\smash{\SetFigFont{12}{14.4}{rm}$f_n$}}}
\put(3826,-2986){\makebox(0,0)[lb]{\smash{\SetFigFont{12}{14.4}{rm}${n\over n+1}$}}}
\end{picture}
 }
\caption[]{}
\label{recu}
\end{figure}

Let $\sigma \in \Sigma(n)$ be an element of the symmetric group on n 
elements, and denote by $T(\sigma)$ the associated tangle in ${\cal 
S}_p(n,n)$.  If $n\in \{0,1,\cdots,p-1\}$, then we define an element 
$f_n$ of ${\cal S}_p{(n,n)}$ as :
\begin{equation}
f_n= \psi_p({1\over n!})\sum_{\sigma \in\Sigma(n)} (-1)^{|\sigma|} T(\sigma).
\end{equation}
In particular $f_0$ is the empty tangle.  This element satisfies the following 
properties~:
\begin{eqnarray}
(f_k\otimes I_{n-k}) \circ f_n =f_n \circ (I_{n-k} \otimes f_k) = f_n,
 \label{ff}\\
T(\sigma) \circ f_n= f_n \circ T(\sigma) =(-1)^{|\sigma|} f_n,\, \sigma 
\in \Sigma_n,\label{perm}\\
g \circ f_n=0 \, (\mbox{resp } f_n\circ g =0) \mbox{ if $ g$ possess}
\nonumber \\
\mbox{ an arc connecting ${\Bbb{R}} \times {0}$
  (resp  ${\Bbb{R}} \times {1}$) with itself},\label{connec}\\
\mbox{ relation in fig \ref{recu} if } n\leq p-2 \label{recur},\\
 trf_n\=p (-1)^n (n+1).\label{trace}
\end{eqnarray}

 The first two properties are clear from the definition of $f_n$; 
for the third property we can suppose, using the behavior of $f_n$ 
under permutations, that $g$ connects the strands $i$ and $i+1$, in which 
case proof is given in figure \ref{nul}, while the first equality is 
obtained from \ref{perm} and the second using fig \ref{writ}.  From 
the definition of $f_n$ it is clear that, for $ n\leq p-2$
\begin{equation}
f_{n+1} = \psi_p({1 \over n+1}) f_n \otimes I -
\psi_p({n\over n+1}) T(t_{n,n+1}) f_n \otimes I
\end{equation}
where $t_{n,n+1}$ is the permutation of $n$ with $n+1$.
Then, using the skein relation and  property \ref{ff} the conclusion 
of \ref{recur} is direct and \ref{trace} is shown by recurrence using 
\ref{recur}.
In particular the relation \ref{trace} implies 
that $f_{p-1} \in \N_p{(n,n)}$.
\begin{figure}[htbp]
\centerline{
\setlength{\unitlength}{0.0125in}
\begingroup\makeatletter\ifx\SetFigFont\undefined
\def\x#1#2#3#4#5#6#7\relax{\def\x{#1#2#3#4#5#6}}%
\expandafter\x\fmtname xxxxxx\relax \def\y{splain}%
\ifx\x\y   
\gdef\SetFigFont#1#2#3{%
  \ifnum #1<17\tiny\else \ifnum #1<20\small\else
  \ifnum #1<24\normalsize\else \ifnum #1<29\large\else
  \ifnum #1<34\Large\else \ifnum #1<41\LARGE\else
     \huge\fi\fi\fi\fi\fi\fi
  \csname #3\endcsname}%
\else
\gdef\SetFigFont#1#2#3{\begingroup
  \count@#1\relax \ifnum 25<\count@\count@25\fi
  \def\x{\endgroup\@setsize\SetFigFont{#2pt}}%
  \expandafter\x
    \csname \romannumeral\the\count@ pt\expandafter\endcsname
    \csname @\romannumeral\the\count@ pt\endcsname
  \csname #3\endcsname}%
\fi
\fi\endgroup
\begin{picture}(141,115)(0,-10)
\drawline(100,85)(110,60)
\drawline(110,85)(100,60)
\drawline(100,85)	(100.349,89.513)
	(100.767,92.829)
	(102.000,97.000)
\drawline(102,97)	(105.000,100.000)
\drawline(105,100)	(108.000,97.000)
\drawline(108,97)	(109.233,92.829)
	(109.651,89.513)
	(110.000,85.000)
\drawline(115,40)(115,60)(75,60)
	(75,40)(115,40)
\drawline(80,100)(80,60)
\drawline(80,40)(80,0)
\drawline(110,40)(110,0)
\drawline(100,40)(100,0)
\put(85,80){\makebox(0,0)[lb]{\smash{{{\SetFigFont{12}{14.4}{rm}...}}}}}
\put(85,20){\makebox(0,0)[lb]{\smash{{{\SetFigFont{12}{14.4}{rm}...}}}}}
\put(93,47){\makebox(0,0)[lb]{\smash{{{\SetFigFont{12}{14.4}{rm}f}}}}}
\drawline(40,40)(40,60)(0,60)
	(0,40)(40,40)
\drawline(5,100)(5,60)
\drawline(5,40)(5,0)
\drawline(35,40)(35,0)
\drawline(25,40)(25,0)
\drawline(25,60)	(24.938,62.616)
	(24.884,65.049)
	(24.800,69.402)
	(24.750,73.132)
	(24.734,76.312)
	(24.750,79.016)
	(24.800,81.316)
	(25.000,85.000)
\drawline(25,85)	(25.211,87.663)
	(25.556,91.014)
	(26.123,94.357)
	(27.000,97.000)
\drawline(27,97)	(30.000,100.000)
\drawline(30,100)	(33.000,97.000)
\drawline(33,97)	(33.877,94.357)
	(34.444,91.014)
	(34.789,87.663)
	(35.000,85.000)
\drawline(35,85)	(35.200,81.316)
	(35.250,79.016)
	(35.266,76.312)
	(35.250,73.132)
	(35.200,69.402)
	(35.116,65.049)
	(35.062,62.616)
	(35.000,60.000)
\put(10,80){\makebox(0,0)[lb]{\smash{{{\SetFigFont{12}{14.4}{rm}...}}}}}
\put(10,20){\makebox(0,0)[lb]{\smash{{{\SetFigFont{12}{14.4}{rm}...}}}}}
\put(18,47){\makebox(0,0)[lb]{\smash{{{\SetFigFont{12}{14.4}{rm}f}}}}}
\put(125,45){\makebox(0,0)[lb]{\smash{{{\SetFigFont{12}{14.4}{rm}= 0}}}}}
\put(62,47){\makebox(0,0)[lb]{\smash{{{\SetFigFont{12}{14.4}{rm}-}}}}}
\put(49,46){\makebox(0,0)[lb]{\smash{{{\SetFigFont{12}{14.4}{rm}=}}}}}
\end{picture}
 }
\caption[]{}
\label{nul}
\end{figure}

The following lemma is satisfied:
\begin{lemma}
Let ${\cal D}_p{(k,l)} \subset {\cal S}_p{(k,l)}$ 
be the set of elements of the form 
\begin{equation}
f = n + \sum_s x_s f_{i_s} y_s. 
\end{equation}
where s runs over a finite set of indices,
${i_s}\in \{0,1,\cdots,p-2\}$, 
$n \in \N_p{(k,l)}$. And $x_s\, (resp.\, y_s)$ is an element of 
${\cal S}_p{(i_s,l)}$ 
(resp.  ${\cal S}_p{(k,i_s)})$.
Then ${\cal S}_p{(k,l)}={\cal D}_p{(k,l)}$.
\end{lemma}
{\bf Proof}

First, note that if $x,y \in {\cal D}_p{(k,l)} $ then $x+y \in {\cal 
D}_p{(k,l)} $ and if $x \in {\cal D}_p{(k,l)} $ $y \in {\cal 
S}_p{(l,n)})$ (resp. $y \in {\cal S}_p{(n,k)}$) then $y\cdot x \in 
{\cal D}_p{(k,n)}$ (resp. $x\cdot y \in {\cal D}_p{(n,l)}$).  So it is 
enough to show this property for the identity tangle $I_n \in {\cal 
S}_p{(n,n)}$.  We will show it by recurrence~: $I_0=f_0, I_1=f_1$, 
suppose that $I_n\in {\cal D}_p{(n,n)}$, i-e $I_n = n + \sum_s x_s 
f_{i_s} y_s $ where $n \in \N_p{(n,n)}$ and $x_s\, (resp.\, y_s)$ is an 
element of ${\cal S}_p{(i_s,n)}$ (resp. ${\cal S}_p{(n,i_s)})$.  Then 
$I_{n+1}= I_n \otimes I_1 = n\otimes I_1 + \sum_s (x_s\otimes I_1) 
(f_{i_s}\otimes I_1) (y_s \otimes I_1) $.  $n\otimes I_1 \in n \in 
\N_p{(n+1,n+1)}$, moreover by relation \ref{recur} $f_{i_s}\otimes I_1 
=f_{i_s +1} + x_{i_s} f_{i_s-1} y_{i_s}$ with $ x_{i_s} \in {\cal 
S}_p{(i_s-1,i_s+1)}$, $y_{i_s}\in {\cal S}_p{(i_s+1,i_s-1)}$.  Thus 
$f_{i_s}\otimes I_1 \in{\cal D}_p{(i_s+1,i_s+1)}$ ($ f_{p-1} \in 
\N_p{(p-1,p-1)}$) , thus $(x_s\otimes I_1) (f_{i_s}\otimes I_1) (y_s 
\otimes I_1) \in {\cal D}_p{(n+1,n+1)}$ and so is $I_{n+1}$.\\
$\Box $

Define $I_p=\{0,1,\cdots,p-2\}$ and $J=\{ {\vec k}=(k_1,\cdots, k_n) 
\, n\in {\Bbb{Z}},k_i \in I_p\}$.  We denote $|{\vec k}| = 
k_1+\cdots+k_n$ and $ f_{\vec k} \in {\cal S}_p({|{\vec k}|,|{\vec 
k}|})$ the element $ f_{k_1} \otimes \cdots \otimes f_{k_n}$.  And we 
define the invariant skein module ${\cal I}{({\vec k},{\vec l})} 
\subset {\cal S}_p{(|{\vec k}|,|{\vec l}|)}$ as follows:
\begin{equation}
{\cal I}{({\vec k},{\vec l})} =\left\{ x \in {\cal S}_p{(|{\vec k}|,|{\vec l}|)} |\,
f_{\vec l}\,x f_{\vec k} =x\right\}.
\end{equation}

We say that the triple $(i,j,k)\in {I_p}^3$ is admissible if and only if
$i+j-k \geq 0 $, $i-j+k \geq 0$, $-i+j+k \geq 0 $ and $i+j+k$ is even.
We say that the triple $(i,j,k)$ is $p$-admissible if and only if it is
admissible and $ {i+j+k} \leq 2(p-2)$.

\begin{lemma}
\label{trival}
The space ${\cal I}{({(i)},{(j,k)})}$, $(i,j,k) \in I_p^3$\\
is a zero dimensional space if 
$(i,j,k)$ is not admissible,\\
is a one dimensional module if 
$(i,j,k)$ is admissible\\
and in this case it is a null submodule (i-e $\subset \N_p{(i,j+k)}$) 
iff $i+j+k \geq 2p-2$.
\end{lemma}
{\bf Proof}
Let $T$ be a simple tangle in $\S(i,j+k)$,and denote 
by $IT= f_j \otimes f_k T f_i $ the associated invariant tangle.  If 
$(i,j,k)$ is not admissible then $T$ possesses an arc connecting $i$ 
(resp.  $j$, $k$) with itself.  Thus the property \ref{connec} implies 
that $IT=0$, so $ {\cal I}{({(i)},{(j,k)})}=0 $. If $(i,j,k)$ is 
admissible there is only one simple diagram $S_i^{(j,k)}$ in 
$\S(i,j+k)$ which do not possess an arc connecting $i$ (resp.  $j$, 
$k$) with itself.  $S_{(i)}^{(j,k)}$ is the simple diagram with $a= 
(-i+j+k)/2$ arcs connecting $j$ with $k$, $b= (i-j+k)/2$ arcs 
connecting $i$ with $k$, $c= (i+j-k)/2$ arcs connecting $i$ with $j$.  
Let $Y_{(i)}^{(j,k)} = f_j \otimes f_k S_i^{(j,k)} f_i $ Then 
$Y_i^{(j,k)}$ is a base element of ${\cal I}{({(i)},{(j,k)})}$ (i-e $ 
{\cal I}{({(i)},{(j,k)})}= K_p Y_{(i)}^{(j,k)}$.  Moreover, if 
$(i,j,k)\in I_p^3 $, then
\begin{equation}
tr(Y_{(j,k)}^{(i)} Y_{(i)}^{(j,k)}) = \psi_p(\Theta(i,j,k)),
\end{equation}
where 
\begin{equation}
\Theta(i,j,k) =(-1)^{{i+j+k \over 2}} 
{({i+j+k\over2}+1)! ({-i+j+k \over 2})!({i-j+k \over 2})!({i+j-k \over 2})! \over (i)! (j)! (k)!}.
\end{equation}
This property is a consequence of \ref{eva}.  Thus 
$\psi_p(\Theta(i,j,k))=0$ if and only if $i+j+k \geq 2(p-1)$, so in 
that case $ Y_{(i)}^{(j,k)} \in N_p{(i,j+k)}$ (resp.  $ 
Y_{(j,k)}^{(i)}\in N_p{(i,j+k)}$).  Moreover we have the following 
identity~:
\begin{equation}
\label{ortho}
Y_{(j,k)}^{(l)} Y_{(i)}^{(j,k)} = 
\psi_{p}({\Theta(i,j,k) \over \Delta_i}) 
\delta_{i,l} f_i,
\end{equation}

if $(i,j,k) \in I_p^3$.
Here $\delta_{i,l}$ denotes the Kronecker symbol and 
$\Delta_i= trf_i = (-1)^i (i+1)$.
Since $Y_{(j,k)}^{(l)} Y_{(i)}^{(j,k)}\in {\cal I}_p{(i,l)}$, by the 
preceding lemma it is equal to $0$ if $i\not=l$ and proportional to 
$f_i$ if $i=l$.  The coefficient of proportionality is calculated by 
taking the trace.\\
$\Box$

\begin{prop}
Let $f \in {\cal I}_p{((k_1,k_2),(l_1,l_2))}$, $k_i,l_i \in I_p$,
then 
\begin{equation}
f =  n + \sum_i \alpha_i Y_{i}^{(l_1,l_2)} Y_{(k_1,k_2)}^{i},
\end{equation}
where $n \in \N_p(k_1+k_2,l_1+l_2)$, the sum
 is over all $i$ such that $(i,l_1,l_2)$ and $(i,k_1,k_2)$
are admissible triples and 
\begin{equation}
\alpha_i = tr(Y_{(l_1,l_2)}^{i} f Y_{i}^{(k_1,k_2)}) 
{\Delta_i  \over \Theta(i,k_1,k_2) \Theta(i,l_1,l_2) }.
\end{equation}
\end{prop}
{\bf Proof}

Let $f \in {\cal I}{((k_1,k_2),(l_1,l_2))}$; using theorem 1 we know 
that $ f =  n + \sum_s x_s f_{i_s} y_s $ where 
 s runs over a finite set of indices,
${i_s}\in I_p $, 
 $n \in \N_p$. And $x_s\, (resp. \,y_s)$ is an element of 
 ${\cal S}_p{(i_s,l_1+l_2)}$  (resp. ${\cal S}_p{(k_1+k_2,i_s)}$).
By composing on the left by $ f_{l_1} \otimes f_{l_2}$,
on the right by $ f_{k_1} \otimes f_{k_2}$ we can assume that 
$x_s \in {\cal I}{((l_1,l_2),i_s)}$, $y_s \in {\cal I}{(i_s,(k_1,k_2))}$.
So, by lemma \ref{trival}, $x_s = \lambda_s Y_{i_s}^{(l_1,l_2)}$ 
if $(i_s,l_1,l_2)$ is $p$-admissible or $x_s \in \N(i_s, l_1 +l_2)$
(resp.  $y_s = \mu_s Y_{(l_1,l_2)}^{i_s}$
if $(i_s,k_1,k_2)$ is $p$-admissible or $y_s \in \N( l_1 +l_2,i_s)$)
\begin{equation}
f =  n + \sum_i \alpha_i Y_{i}^{(l_1,l_2)} Y_{(k_1,k_2)}^{i}  .
\end{equation}
The sum is over all $i$ such that $(i,l_1,l_2)$ and $(i,k_1,k_2)$
are $p$-admissible triples and $\alpha_i \in K_p$.
Moreover, multiplying both side by $Y_{(l_1,l_2)}^{i}$ on the left 
and $Y_{i}^{(k_1,k_2)}$ on the right and  using \ref{ortho}, 
we get the desired result for  $\alpha_i$.\\
$\Box$

In particular, 
\begin{equation}
\label{identity}
f_i \otimes f_j = n+ \sum_{k \atop (i,j,k) p-admissible}
\psi_p({\Delta_k \over \Theta(i,j,k)})  Y_{k}^{(i,j)} Y_{(i,j)}^{k}.
\end{equation}

\section{Spin Networks, Chord diagrams and weight system} 
\begin{figure}[htbp]
\centerline{
\begin{picture}(0,0)%
\epsfig{file=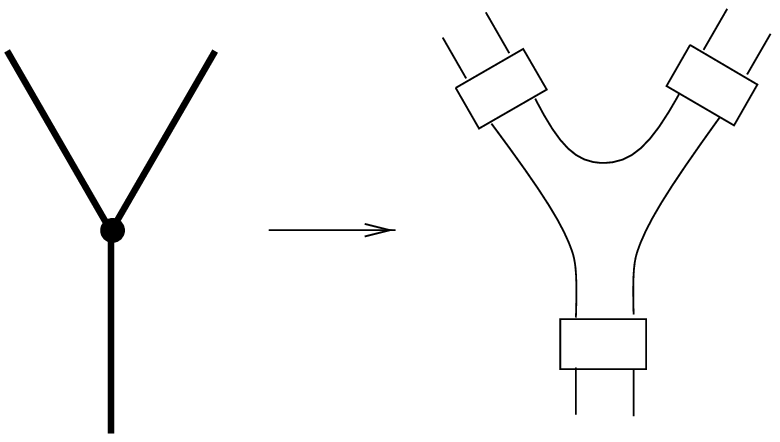}%
\end{picture}%
\setlength{\unitlength}{0.00083300in}%
\begingroup\makeatletter\ifx\SetFigFont\undefined
\def\x#1#2#3#4#5#6#7\relax{\def\x{#1#2#3#4#5#6}}%
\expandafter\x\fmtname xxxxxx\relax \def\y{splain}%
\ifx\x\y   
\gdef\SetFigFont#1#2#3{%
  \ifnum #1<17\tiny\else \ifnum #1<20\small\else
  \ifnum #1<24\normalsize\else \ifnum #1<29\large\else
  \ifnum #1<34\Large\else \ifnum #1<41\LARGE\else
     \huge\fi\fi\fi\fi\fi\fi
  \csname #3\endcsname}%
\else
\gdef\SetFigFont#1#2#3{\begingroup
  \count@#1\relax \ifnum 25<\count@\count@25\fi
  \def\x{\endgroup\@setsize\SetFigFont{#2pt}}%
  \expandafter\x
    \csname \romannumeral\the\count@ pt\expandafter\endcsname
    \csname @\romannumeral\the\count@ pt\endcsname
  \csname #3\endcsname}%
\fi
\fi\endgroup
\begin{picture}(3709,2082)(669,-5269)
\put(751,-3833){\makebox(0,0)[lb]{\smash{\SetFigFont{12}{14.4}{rm}j}}}
\put(1583,-3848){\makebox(0,0)[lb]{\smash{\SetFigFont{12}{14.4}{rm}k}}}
\put(1028,-4898){\makebox(0,0)[lb]{\smash{\SetFigFont{12}{14.4}{rm}i}}}
\put(3424,-4854){\makebox(0,0)[lb]{\smash{\SetFigFont{12}{14.4}{rm}$f_i$}}}
\put(3926,-3517){\makebox(0,0)[lb]{\smash{
\SetFigFont{12}{14.4}{rm}$f_k$
}}}
\put(2946,-3704){\makebox(0,0)[lb]{\smash{
\SetFigFont{12}{14.4}{rm}$f_j$
}}}
\end{picture}
 }
\caption[]{}
\label{trivalentvert}
\end{figure}

A $p$-spin 
network is a trivalent graph $\Gamma$ equipped with a cyclic 
orientation at each vertex and with an $p$-admissible coloring,
which is a mapping  from the edges of $\Gamma$ to $I_p$, such that 
every triple $(i,j,k)\in I_p^3$ surrounding a trivalent vertex is 
$p$-admissible.  If $\Gamma$ is a $p$-spin network we can associate  
a closed invariant tangle, denoted $e_p(\Gamma) \in 
\I_p(0,0)={{\Bbb{Z}}/ p{\Bbb{Z}}}$, to any immersion of $\Gamma$ into 
${\Bbb{R}} \times [0,1] $.  The correspondence is the following : To 
each edge colored by $i\in I_p$ of the immersed spin network we 
associate the tangle $f_i \in \I(i,i)$, and to each trivalent vertex 
of the immersed spin network colored by the admissible triple 
$(i,j,k)$, with a positive orientation we associate the tangle 
$Y(i,j,k)$ (see figure \ref{trivalentvert}).
Moreover, the relations of lemma \ref{isotopy}
imply that this correspondence does not 
depend on the particular chosen immersion of the spin network.

Untill now we have worked in characteristic $p$, but of course it is 
much more standard to work in characteristic $0$ \cite{Kaufflin, 
turaevb}, in which case we just have to replace in what we said $ {{\Bbb{Z}}/ 
p{\Bbb{Z}}}$ by ${\Bbb{Q}}$, $\psi_p$ by the identity map, $I_p$ by 
$\Bbb{N}$ and so on.  Let $e$ be the standard spin network evaluation 
in characteristic $0$, then~:
\begin{prop}
\label{caracp}
If $\Gamma$ is a $p$-spin network 
\begin{equation}
    \label{eva}
    e_p(\Gamma)=\psi_p(e(\Gamma))
\end{equation}
\end{prop}
{\bf Proof}\\
If $\Gamma$ is a $p$-spin network,
using the definition of $f_i$ we can expand $e(\Gamma)$ into a 
 finite sum $\sum_i a_i T^{(i)}$, where $ T^{(i)}$ are closed tangles and 
$a_i \in {\Bbb{Q}}_p$. If we evaluate $e_p(\Gamma)$ we get the same result
but with $a_i$ replaced by $\psi_p(a_i)$.
 Each closed tangle is evaluated as an integer 
 using \ref{circ}, \ref{skei}, 
these relations don't depend on the characteristic thus we get the proposition
\ref{caracp}.\\
$\Box$

\begin{figure}[htbp]
\centerline{
\begin{picture}(0,0)%
\epsfig{file=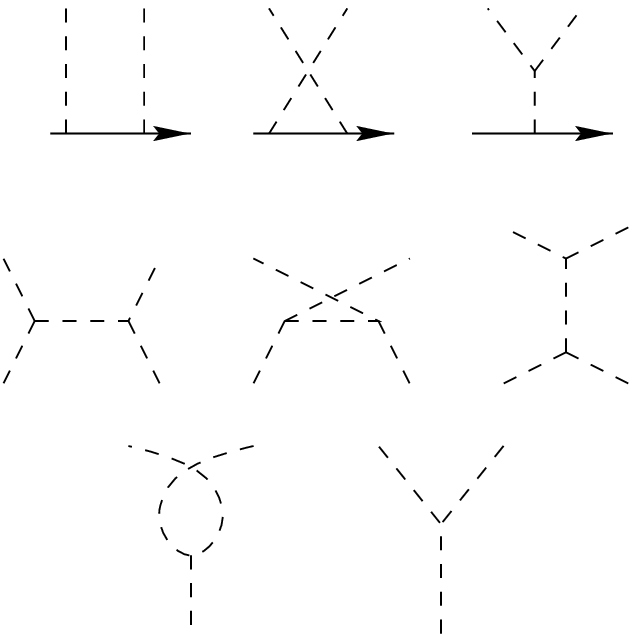}%
\end{picture}%
\setlength{\unitlength}{0.00083300in}%
\begingroup\makeatletter\ifx\SetFigFont\undefined
\def\x#1#2#3#4#5#6#7\relax{\def\x{#1#2#3#4#5#6}}%
\expandafter\x\fmtname xxxxxx\relax \def\y{splain}%
\ifx\x\y   
\gdef\SetFigFont#1#2#3{%
  \ifnum #1<17\tiny\else \ifnum #1<20\small\else
  \ifnum #1<24\normalsize\else \ifnum #1<29\large\else
  \ifnum #1<34\Large\else \ifnum #1<41\LARGE\else
     \huge\fi\fi\fi\fi\fi\fi
  \csname #3\endcsname}%
\else
\gdef\SetFigFont#1#2#3{\begingroup
  \count@#1\relax \ifnum 25<\count@\count@25\fi
  \def\x{\endgroup\@setsize\SetFigFont{#2pt}}%
  \expandafter\x
    \csname \romannumeral\the\count@ pt\expandafter\endcsname
    \csname @\romannumeral\the\count@ pt\endcsname
  \csname #3\endcsname}%
\fi
\fi\endgroup
\begin{picture}(3024,3024)(589,-3073)
\put(2566,-716){\makebox(0,0)[lb]{\smash{\SetFigFont{12}{14.4}{rm}$=$}}}
\put(1523,-716){\makebox(0,0)[lb]{\smash{\SetFigFont{12}{14.4}{rm}$-$}}}
\put(1448,-1616){\makebox(0,0)[lb]{\smash{\SetFigFont{12}{14.4}{rm}$-$}}}
\put(2566,-1616){\makebox(0,0)[lb]{\smash{\SetFigFont{12}{14.4}{rm}$=$}}}
\put(1801,-2611){\makebox(0,0)[lb]{\smash{\SetFigFont{12}{14.4}{rm}$= -$}}}
\end{picture}
}
\caption[]{}
\label{ihx}
\end{figure}

Let $X$ be a one dimensional oriented compact manifold without boundary.
A chord diagram (usually refered as Chinese character chord diagram)
is the union $D =\bar{D}\cup X$ where $\bar{D}$ is a graph with univalent
 and trivalent vertices, together with a cyclic orientation of trivalent
 vertices  such that univalent vertices lie in $X$.  Trivalent 
 vertices are referred to as internal vertices, and the degree of $D$, 
 denoted $d^{\circ}(D)$, is half the number of vertices of the graph $ 
 \bar{D}$.  Let $\tilde {\cal A}_n$ the $\Bbb{Z}$ module freely 
 generated by chord diagrams of degree $n$.  We define the ${\Bbb{Z}}$ 
 module of chord diagrams of degree n, denoted ${\cal A}_n$, as being the 
 quotient of $\tilde {\cal A}_n$ by the relations (STU, IHX, AS) shown 
 in figure \ref{ihx}.

\begin{figure}[htbp]
\centerline{
\begin{picture}(0,0)%
\epsfig{file=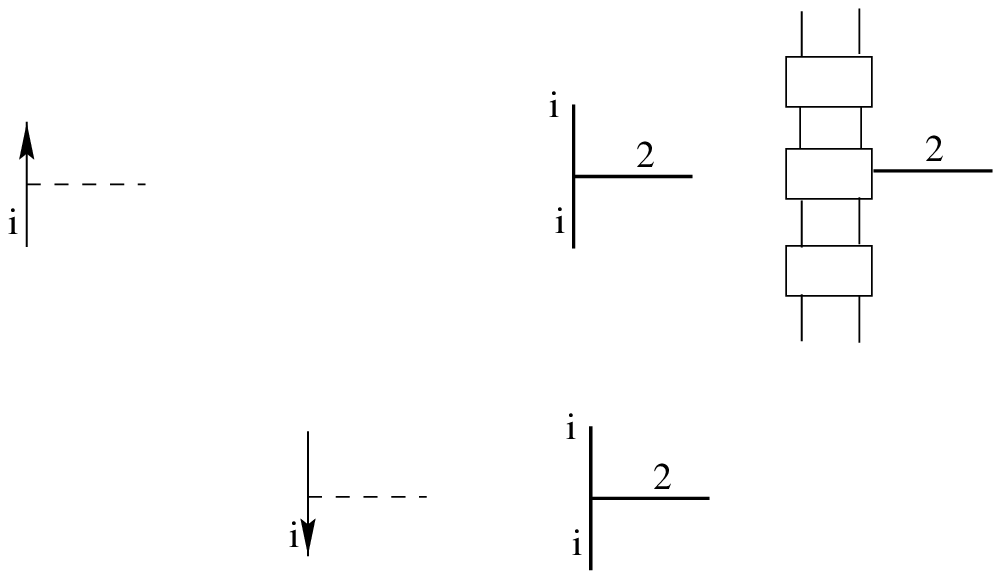}%
\end{picture}%
\setlength{\unitlength}{0.00083300in}%
\begingroup\makeatletter\ifx\SetFigFont\undefined
\def\x#1#2#3#4#5#6#7\relax{\def\x{#1#2#3#4#5#6}}%
\expandafter\x\fmtname xxxxxx\relax \def\y{splain}%
\ifx\x\y   
\gdef\SetFigFont#1#2#3{%
  \ifnum #1<17\tiny\else \ifnum #1<20\small\else
  \ifnum #1<24\normalsize\else \ifnum #1<29\large\else
  \ifnum #1<34\Large\else \ifnum #1<41\LARGE\else
     \huge\fi\fi\fi\fi\fi\fi
  \csname #3\endcsname}%
\else
\gdef\SetFigFont#1#2#3{\begingroup
  \count@#1\relax \ifnum 25<\count@\count@25\fi
  \def\x{\endgroup\@setsize\SetFigFont{#2pt}}%
  \expandafter\x
    \csname \romannumeral\the\count@ pt\expandafter\endcsname
    \csname @\romannumeral\the\count@ pt\endcsname
  \csname #3\endcsname}%
\fi
\fi\endgroup
\begin{picture}(5137,2730)(376,-6960)
\put(376,-5139){\makebox(0,0)[lb]{\smash{\SetFigFont{12}{14.4}{rm}$D$ $=$}}}
\put(4606,-5096){\makebox(0,0)[lb]{\smash{\SetFigFont{12}{14.4}{rm}$\Delta$}}}
\put(4592,-4654){\makebox(0,0)[lb]{\smash{\SetFigFont{12}{14.4}{rm}$f_i$}}}
\put(4585,-5547){\makebox(0,0)[lb]{\smash{\SetFigFont{12}{14.4}{rm}$f_i$}}}
\put(4201,-5086){\makebox(0,0)[lb]{\smash{\SetFigFont{12}{14.4}{rm}$=$}}}
\put(1868,-5093){\makebox(0,0)[lb]{\smash{\SetFigFont{12}{14.4}{rm}$\omega^{T}(D) =$}}}
\put(2709,-5085){\makebox(0,0)[lb]{\smash{\SetFigFont{12}{14.4}{rm} $i \cdot e($}}}
\put(4066,-5096){\makebox(0,0)[lb]{\smash{\SetFigFont{12}{14.4}{rm}$)$}}}
\put(2791,-6638){\makebox(0,0)[lb]{\smash{\SetFigFont{12}{14.4}{rm}$)$}}}
\put(2948,-6653){\makebox(0,0)[lb]{\smash{\SetFigFont{12}{14.4}{rm}$=$ $-i \cdot$}}}
\put(1748,-6635){\makebox(0,0)[lb]{\smash{\SetFigFont{12}{14.4}{rm}$\bar{\omega}^{T}($}}}
\end{picture}
 }
\caption[]{}
\label{rule}
\end{figure}

\begin{figure}[htbp]
\centerline{
\begin{picture}(0,0)%
\epsfig{file=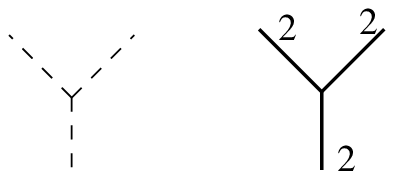}%
\end{picture}%
\setlength{\unitlength}{0.00083300in}%
\begingroup\makeatletter\ifx\SetFigFont\undefined
\def\x#1#2#3#4#5#6#7\relax{\def\x{#1#2#3#4#5#6}}%
\expandafter\x\fmtname xxxxxx\relax \def\y{splain}%
\ifx\x\y   
\gdef\SetFigFont#1#2#3{%
  \ifnum #1<17\tiny\else \ifnum #1<20\small\else
  \ifnum #1<24\normalsize\else \ifnum #1<29\large\else
  \ifnum #1<34\Large\else \ifnum #1<41\LARGE\else
     \huge\fi\fi\fi\fi\fi\fi
  \csname #3\endcsname}%
\else
\gdef\SetFigFont#1#2#3{\begingroup
  \count@#1\relax \ifnum 25<\count@\count@25\fi
  \def\x{\endgroup\@setsize\SetFigFont{#2pt}}%
  \expandafter\x
    \csname \romannumeral\the\count@ pt\expandafter\endcsname
    \csname @\romannumeral\the\count@ pt\endcsname
  \csname #3\endcsname}%
\fi
\fi\endgroup
\begin{picture}(1995,803)(1628,-7048)
\put(2341,-6713){\makebox(0,0)[lb]{\smash{\SetFigFont{12}{14.4}{rm}$)$}}}
\put(2565,-6713){\makebox(0,0)[lb]{\smash{\SetFigFont{12}{14.4}{rm}$=$ $2 \cdot$}}}
\put(1628,-6710){\makebox(0,0)[lb]{\smash{\SetFigFont{12}{14.4}{rm}$\bar{\omega}^{T}($}}}
\end{picture}
 }
\caption[]{}
\label{weightint}
\end{figure}

\begin{figure}[htbp]
\centerline{
\begin{picture}(0,0)%
\epsfig{file=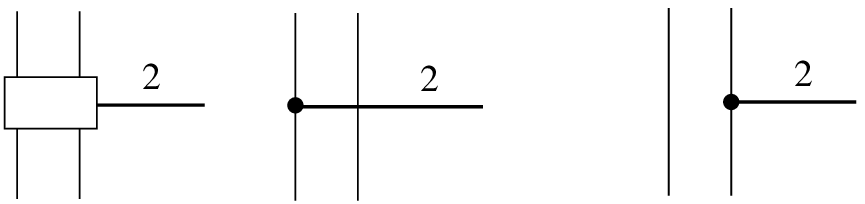}%
\end{picture}%
\setlength{\unitlength}{0.00083300in}%
\begingroup\makeatletter\ifx\SetFigFont\undefined
\def\x#1#2#3#4#5#6#7\relax{\def\x{#1#2#3#4#5#6}}%
\expandafter\x\fmtname xxxxxx\relax \def\y{splain}%
\ifx\x\y   
\gdef\SetFigFont#1#2#3{%
  \ifnum #1<17\tiny\else \ifnum #1<20\small\else
  \ifnum #1<24\normalsize\else \ifnum #1<29\large\else
  \ifnum #1<34\Large\else \ifnum #1<41\LARGE\else
     \huge\fi\fi\fi\fi\fi\fi
  \csname #3\endcsname}%
\else
\gdef\SetFigFont#1#2#3{\begingroup
  \count@#1\relax \ifnum 25<\count@\count@25\fi
  \def\x{\endgroup\@setsize\SetFigFont{#2pt}}%
  \expandafter\x
    \csname \romannumeral\the\count@ pt\expandafter\endcsname
    \csname @\romannumeral\the\count@ pt\endcsname
  \csname #3\endcsname}%
\fi
\fi\endgroup
\begin{picture}(4122,947)(461,-373)
\put(1471, 36){\makebox(0,0)[lb]{\smash{\SetFigFont{12}{14.4}{rm}$=$}}}
\put(548, 29){\makebox(0,0)[lb]{\smash{\SetFigFont{12}{14.4}{rm}$\Delta$}}}
\put(2844, 37){\makebox(0,0)[lb]{\smash{\SetFigFont{12}{14.4}{rm}$+$}}}
\put(3368, 38){\makebox(0,0)[lb]{\smash{\SetFigFont{12}{14.4}{rm}$+$}}}
\put(3676,328){\makebox(0,0)[lb]{\smash{\SetFigFont{12}{14.4}{rm}$\cdots$}}}
\put(3098, 36){\makebox(0,0)[lb]{\smash{\SetFigFont{12}{14.4}{rm}$\cdots$}}}
\put(1891,336){\makebox(0,0)[lb]{\smash{\SetFigFont{12}{14.4}{rm}$\cdots$}}}
\put(571,329){\makebox(0,0)[lb]{\smash{\SetFigFont{12}{14.4}{rm}$\cdots$}}}
\end{picture}
 }
\caption[]{}
\label{coproduit}
\end{figure}

We denote by $\vec{\lambda}$ a $p$-coloring of $X$ i-e a mapping from 
the set of connected components of $X$ to $ I_p=\{1,\cdots, p-2\}$.  
And we 
define a weight system $ \bar{\omega}^{T}_{\vec{\lambda}}$ which 
associates a spin network to a chord diagram.  The rules defining $ 
\omega^{T}$ are given in the figures (\ref{rule},\ref{weightint}) 
where $\Delta$ is the coproduct, as shown in figure \ref{coproduit}.

\begin{figure}[htbp]
\centerline{
\begin{picture}(0,0)%
\epsfig{file=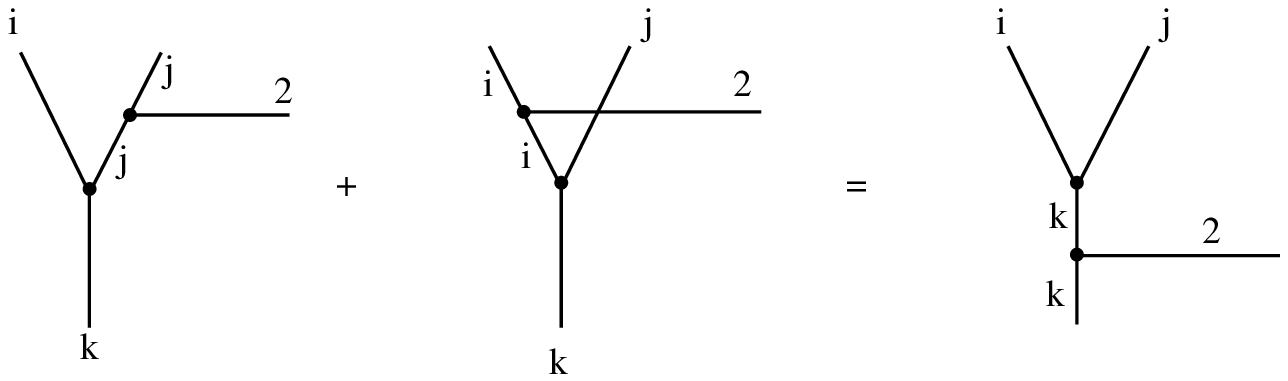}%
\end{picture}%
\setlength{\unitlength}{0.00083300in}%
\begingroup\makeatletter\ifx\SetFigFont\undefined
\def\x#1#2#3#4#5#6#7\relax{\def\x{#1#2#3#4#5#6}}%
\expandafter\x\fmtname xxxxxx\relax \def\y{splain}%
\ifx\x\y   
\gdef\SetFigFont#1#2#3{%
  \ifnum #1<17\tiny\else \ifnum #1<20\small\else
  \ifnum #1<24\normalsize\else \ifnum #1<29\large\else
  \ifnum #1<34\Large\else \ifnum #1<41\LARGE\else
     \huge\fi\fi\fi\fi\fi\fi
  \csname #3\endcsname}%
\else
\gdef\SetFigFont#1#2#3{\begingroup
  \count@#1\relax \ifnum 25<\count@\count@25\fi
  \def\x{\endgroup\@setsize\SetFigFont{#2pt}}%
  \expandafter\x
    \csname \romannumeral\the\count@ pt\expandafter\endcsname
    \csname @\romannumeral\the\count@ pt\endcsname
  \csname #3\endcsname}%
\fi
\fi\endgroup
\begin{picture}(6487,1806)(2266,-4273)
\put(2266,-3361){\makebox(0,0)[lb]{\smash{\SetFigFont{12}{14.4}{rm}$j \cdot e($}}}
\put(4411,-3376){\makebox(0,0)[lb]{\smash{\SetFigFont{12}{14.4}{rm}$i \cdot e ($}}}
\put(6961,-3361){\makebox(0,0)[lb]{\smash{\SetFigFont{12}{14.4}{rm}$k \cdot e ($}}}
\put(3586,-3370){\makebox(0,0)[lb]{\smash{\SetFigFont{12}{14.4}{rm}$)$}}}
\put(5844,-3362){\makebox(0,0)[lb]{\smash{\SetFigFont{12}{14.4}{rm}$)$}}}
\put(8462,-3317){\makebox(0,0)[lb]{\smash{\SetFigFont{12}{14.4}{rm}$)$}}}
\end{picture}
 }
\caption[]{}
\label{coprod}
\end{figure}

\begin{figure}[htbp]
\centerline{
\begin{picture}(0,0)%
\epsfig{file=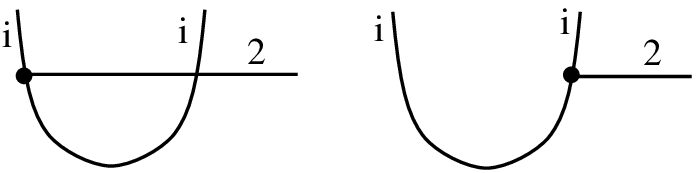}%
\end{picture}%
\setlength{\unitlength}{0.00083300in}%
\begingroup\makeatletter\ifx\SetFigFont\undefined
\def\x#1#2#3#4#5#6#7\relax{\def\x{#1#2#3#4#5#6}}%
\expandafter\x\fmtname xxxxxx\relax \def\y{splain}%
\ifx\x\y   
\gdef\SetFigFont#1#2#3{%
  \ifnum #1<17\tiny\else \ifnum #1<20\small\else
  \ifnum #1<24\normalsize\else \ifnum #1<29\large\else
  \ifnum #1<34\Large\else \ifnum #1<41\LARGE\else
     \huge\fi\fi\fi\fi\fi\fi
  \csname #3\endcsname}%
\else
\gdef\SetFigFont#1#2#3{\begingroup
  \count@#1\relax \ifnum 25<\count@\count@25\fi
  \def\x{\endgroup\@setsize\SetFigFont{#2pt}}%
  \expandafter\x
    \csname \romannumeral\the\count@ pt\expandafter\endcsname
    \csname @\romannumeral\the\count@ pt\endcsname
  \csname #3\endcsname}%
\fi
\fi\endgroup
\begin{picture}(3452,790)(3526,-6536)
\put(6978,-6146){\makebox(0,0)[lb]{\smash{\SetFigFont{12}{14.4}{rm}$=$ $0$}}}
\put(5022,-6108){\makebox(0,0)[lb]{\smash{\SetFigFont{12}{14.4}{rm}$+$}}}
\end{picture}
 }
\caption[]{}
\label{cup}
\end{figure}

\begin{figure}[htbp]
\centerline{
\begin{picture}(0,0)%
\epsfig{file=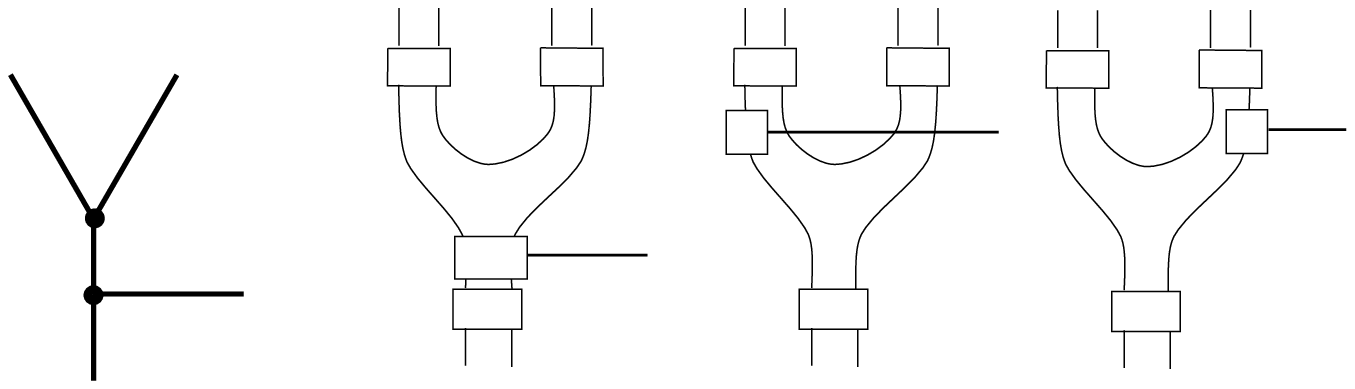}%
\end{picture}%
\setlength{\unitlength}{0.00066700in}%
\begingroup\makeatletter\ifx\SetFigFont\undefined
\def\x#1#2#3#4#5#6#7\relax{\def\x{#1#2#3#4#5#6}}%
\expandafter\x\fmtname xxxxxx\relax \def\y{splain}%
\ifx\x\y   
\gdef\SetFigFont#1#2#3{%
  \ifnum #1<17\tiny\else \ifnum #1<20\small\else
  \ifnum #1<24\normalsize\else \ifnum #1<29\large\else
  \ifnum #1<34\Large\else \ifnum #1<41\LARGE\else
     \huge\fi\fi\fi\fi\fi\fi
  \csname #3\endcsname}%
\else
\gdef\SetFigFont#1#2#3{\begingroup
  \count@#1\relax \ifnum 25<\count@\count@25\fi
  \def\x{\endgroup\@setsize\SetFigFont{#2pt}}%
  \expandafter\x
    \csname \romannumeral\the\count@ pt\expandafter\endcsname
    \csname @\romannumeral\the\count@ pt\endcsname
  \csname #3\endcsname}%
\fi
\fi\endgroup
\begin{picture}(8235,2281)(503,-5269)
\put(3031,-3409){\makebox(0,0)[lb]{\smash{\SetFigFont{10}{12.0}{rm}$f_j$}}}
\put(3931,-3415){\makebox(0,0)[lb]{\smash{\SetFigFont{10}{12.0}{rm}$f_k$}}}
\put(3424,-4854){\makebox(0,0)[lb]{\smash{\SetFigFont{10}{12.0}{rm}$f_i$}}}
\put(5108,-3409){\makebox(0,0)[lb]{\smash{\SetFigFont{10}{12.0}{rm}$f_j$}}}
\put(6008,-3415){\makebox(0,0)[lb]{\smash{\SetFigFont{10}{12.0}{rm}$f_k$}}}
\put(5501,-4854){\makebox(0,0)[lb]{\smash{\SetFigFont{10}{12.0}{rm}$f_i$}}}
\put(6983,-3422){\makebox(0,0)[lb]{\smash{\SetFigFont{10}{12.0}{rm}$f_j$}}}
\put(7883,-3428){\makebox(0,0)[lb]{\smash{\SetFigFont{10}{12.0}{rm}$f_k$}}}
\put(7376,-4867){\makebox(0,0)[lb]{\smash{\SetFigFont{10}{12.0}{rm}$f_i$}}}
\put(751,-3833){\makebox(0,0)[lb]{\smash{\SetFigFont{10}{12.0}{rm}j}}}
\put(1583,-3848){\makebox(0,0)[lb]{\smash{\SetFigFont{10}{12.0}{rm}k}}}
\put(1028,-4898){\makebox(0,0)[lb]{\smash{\SetFigFont{10}{12.0}{rm}i}}}
\put(1546,-4611){\makebox(0,0)[lb]{\smash{\SetFigFont{10}{12.0}{rm}
2
}}}
\put(3443,-4558){\makebox(0,0)[lb]{\smash{\SetFigFont{10}{12.0}{rm}
$\Delta$
}}}
\put(4321,-4423){\makebox(0,0)[lb]{\smash{\SetFigFont{10}{12.0}{rm}
2
}}}
\put(503,-4551){\makebox(0,0)[lb]{\smash{\SetFigFont{10}{12.0}{rm}
$i \cdot$
}}}
\put(2258,-4311){\makebox(0,0)[lb]{\smash{\SetFigFont{10}{12.0}{rm}
$=$
}}}
\put(6504,-3685){\makebox(0,0)[lb]{\smash{\SetFigFont{10}{12.0}{rm}
2
}}}
\put(8582,-3685){\makebox(0,0)[lb]{\smash{\SetFigFont{10}{12.0}{rm}
2
}}}
\put(6481,-4270){\makebox(0,0)[lb]{\smash{\SetFigFont{10}{12.0}{rm}
$+$
}}}
\put(4719,-4293){\makebox(0,0)[lb]{\smash{\SetFigFont{10}{12.0}{rm}
$=$
}}}
\put(5004,-3797){\makebox(0,0)[lb]{\smash{\SetFigFont{10}{12.0}{rm}
$\Delta$
}}}
\put(7989,-3791){\makebox(0,0)[lb]{\smash{\SetFigFont{10}{12.0}{rm}
$\Delta$
}}}
\end{picture}
 }
\caption[]{}
\label{distribution}
\end{figure}

Using the tangle evaluation $e$ of spin networks 
we define the linear form :
\begin{equation}
{\omega}^{T}_{\vec{\lambda}}= e\circ \bar{\omega}^{T}_{\vec{\lambda}}
 :\tilde {\cal A} \longrightarrow {\Bbb{Q}}_p
\end{equation}
We have the following proposition~:
\begin{prop}
${\omega}^{T}_{\vec{\lambda}}$ is a weight system 
i-e it defines a linear form on $\cal A$.
\end{prop}
{\bf Proof}\\
This  is a direct consequence of the relation presented in figure \ref{coprod}
which is itself consequence of the relations of figures \ref{cup} and 
\ref{distribution}.

It is well known that given a Lie algebra $\cal G$ and  an 
invariant  scalar product on $\cal G$ we can associate a weight system
 to any coloring of Wilson lines by representations of $\cal G$.
We denote by $(e,h,f)$ the basis elements of $sl(2)$ which satisfy the 
following commutation relations:
\begin{eqnarray}
&[h,e]=2e, \, [h,f] =-2f  \\
&[e,f]=h.
\end{eqnarray}
The quadratic Casimir is given by $ C =e\otimes f +f\otimes e + 
{1\over 2} h\otimes h$.  We denote by $V_{\lambda}$ the irreducible 
highest weight $sl(2)$ module of weight $\lambda$ given in a basis 
$(v_i),i\in \{0,1,\cdots, \lambda\}$ by:
\begin{eqnarray}
\label{rep}
&e v_i = (\lambda -i+1) v_{i-1} \\
&f v_i = (i+1) v_{i+1}  \\
&h v_i = (\lambda -2i) v_{i} \\
&e v_0 =0, f v_\lambda = 0.
\label{sl2}
\end{eqnarray}
The value of the quadratic Casimir in the representation $V_{\lambda}$ 
is given by $C_\lambda ={\lambda(\lambda +2) \over 2}$. 
Let $\vec{\lambda}=(\lambda_1, \cdots, \lambda_n)\in \Bbb{N}^n$ and 
define $\omega^{sl(2)}_{\vec{\lambda}}$ the $sl(2)$ weight system 
associated with the Lie algebra $sl(2)$ and with a coloring of of the 
components of $(S^1)^L$ by finite dimensional irreducible 
representations $V_{\lambda_i}$ of $sl(2)$ and with the normalization 
of the quadratic casimir given by the trace in the fundamental 
representation.  Then~:
\begin{lemma}
\begin{equation}
\omega^{{sl(2)}}_{\vec{\lambda}}(D) = 
(-1)^{(\lambda_1+\cdots \lambda_n)}(-1)^{deg(D)} 
\omega^{\mathrm T}_{\vec{\lambda}}(D).
\end{equation}
\end{lemma}
{\bf Proof}\\
 Denote by ${\cal T} (k,l) $ the ${\Bbb{Z}}$ linear span of all $(k,l)$ 
 tangles and by $Hom_{sl(2)}(V_1^{\otimes k} ,V_1^{\otimes l})$ the space 
 of homomorphism from $V_1^{\otimes k}$ to $V_1^{\otimes l}$ which 
 commute with the action of $sl(2)$, where the action of $sl(2)$ on 
 $V_1^{\otimes k}$ is the usual diagonal action and the action of 
 $sl(2)$ on $V_1$ is given by \ref{rep}.  We define a linear map $\tau 
 : {\cal T} (k,l) \rightarrow Hom_{sl(2)}(V_1^{\otimes k} ,V_1^{\otimes 
 l})$ as follows~: Let $I_1$ be the $(1,1)$ tangle consisting of one 
 vertical line, $X$ be the $(2,2)$ tangle consisting in a simple 
 crossing, $\cap$ (resp.  $\cup$) be the simplest $(2,0)$ (resp.  
 $(0,2)$) tangle consisting of one line and no crossing, and define
\begin{eqnarray}
\tau(I)&=&Id_{V_1} \\
\tau(X)&=& - P:v\otimes v^{\prime} \rightarrow - v^{\prime}\otimes v \\
\tau(\cap) &:& v\otimes v^{\prime} \rightarrow (v|\epsilon(v^\prime)) \\
\tau(\cup) &:& \alpha \rightarrow  \alpha(v_0\otimes v_1 -v_1\otimes v_0)
\end{eqnarray}
where $(v_0, v_1)$ is the basis \ref{rep} of $V_1$, 
$(v_i|v_j)=\delta_{i,j}$ and $\epsilon(v_0)=v_1, \epsilon(v_1)=-v_0$.
Conjugation by $\epsilon$ of an element  $x \in sl(2)$ 
is the automorphism $*$, $e^*=-f, f^*=-e,h^*=-h$.
We extend $\tau$ to all $(k,l)$ tangles by asking that
$\tau(T\cdot T^{\prime})=\tau (T) \tau(T^{\prime})$ and 
$\tau(T\otimes T^{\prime})=\tau (T) \otimes \tau(T^{\prime})$.
By a direct computation we verify that 
$\tau(O)=-2$, where $O$ denote the circle and 
$\tau(I\otimes I) +\tau(X) + \tau(\cup \cdot \cap)=0$,
so that $\tau$ is well defined on the skein module $\S$.  Let $C_1 \in 
Hom(V_1\otimes V_1)$ be action of the quadratic Casimir on $V_1\otimes 
V_1$, then $C_1 ={1\over 2} \tau (\cup \cdot \cap) -{1\over 2} 
\tau(X)$.  Now if $t$ is the chord diagram consisting of one 
horizontal chord connecting 2 vertical strands then, by the definition 
of the $sl(2)$-tangle weight system, we have that $\omega^{sl(2)}(t) =- 
\tau(\omega^{T}(t)) =C_1$.  Moreover $\tau(f_i)$ corresponds to the 
symetrisation operation and thus  to the projection operator $P_i$ 
from $V_1^{\otimes i}$ to the irreducible representation $V_i$.  If 
$f\in {\cal I}(i,i)$ is an invariant tangle then using the definitions 
of the traces and of $\cap, \cup$ we see that $tr(f)=(-1)^i 
tr_{V_i}(\tau{f})$.  Moreover, $\omega_{i,j}^{sl(2)}(t) 
=c_{i,j}=P_i\otimes P_j \Delta^{(i)}\otimes \Delta^{(j)} C P_i\otimes 
P_j $ which is equal to $-\tau(\omega^T_{i,j}(t))$.\\
$\Box$

\section{$sl(2)-so(3)$ Kirby Weight systems}

 Let $D$ be a chord diagram which is built on an oriented 
 one-dimensional compact manifold $X$, and denote $\hat C$ the chord 
 diagram obtained from $C$ by reversing the orientation of one 
 connected component $X_i$ of $X$.  A weight system $\omega$ is said 
 independent of the orientation if
\begin{equation}
\omega (D) = (-1)^{n_i(D)} \omega (\hat D)
\end{equation}
 where $n_i(D)$ is the number of vertices  lying on $X_i$.
 
 \begin{defi}
 A weight system is a Kirby weight system if it 
 is independent of the orientation 
and if it takes the same value on any two chord diagrams 
 which are related as in figure \ref{Kirby}.
 \end{defi}
 
\begin{figure}[htbp]
\centerline{
\setlength{\unitlength}{0.0125in}
\begingroup\makeatletter\ifx\SetFigFont\undefined
\def\x#1#2#3#4#5#6#7\relax{\def\x{#1#2#3#4#5#6}}%
\expandafter\x\fmtname xxxxxx\relax \def\y{splain}%
\ifx\x\y   
\gdef\SetFigFont#1#2#3{%
  \ifnum #1<17\tiny\else \ifnum #1<20\small\else
  \ifnum #1<24\normalsize\else \ifnum #1<29\large\else
  \ifnum #1<34\Large\else \ifnum #1<41\LARGE\else
     \huge\fi\fi\fi\fi\fi\fi
  \csname #3\endcsname}%
\else
\gdef\SetFigFont#1#2#3{\begingroup
  \count@#1\relax \ifnum 25<\count@\count@25\fi
  \def\x{\endgroup\@setsize\SetFigFont{#2pt}}%
  \expandafter\x
    \csname \romannumeral\the\count@ pt\expandafter\endcsname
    \csname @\romannumeral\the\count@ pt\endcsname
  \csname #3\endcsname}%
\fi
\fi\endgroup
\begin{picture}(257,95)(0,-10)
\put(84,47){\makebox(0,0)[lb]{\smash{{{\SetFigFont{20}{24.0}{rm}.}}}}}
\drawline(40.000,66.000)(34.928,65.500)(30.050,64.021)
	(25.555,61.618)(21.615,58.385)(18.382,54.445)
	(15.979,49.950)(14.500,45.072)(14.000,40.000)
	(14.500,34.928)(15.979,30.050)(18.382,25.555)
	(21.615,21.615)(25.555,18.382)(30.050,15.979)
	(34.928,14.500)(40.000,14.000)
\drawline(40.000,14.000)(45.212,14.328)(50.254,15.690)
	(54.922,18.032)(59.029,21.259)(62.408,25.240)
	(64.924,29.817)(66.475,34.804)(67.000,40.000)
	(66.475,45.196)(64.924,50.183)(62.408,54.760)
	(59.029,58.741)(54.922,61.968)(50.254,64.310)
	(45.212,65.672)(40.000,66.000)
\drawline(190.000,66.000)(184.928,65.500)(180.050,64.021)
	(175.555,61.618)(171.615,58.385)(168.382,54.445)
	(165.979,49.950)(164.500,45.072)(164.000,40.000)
	(164.500,34.928)(165.979,30.050)(168.382,25.555)
	(171.615,21.615)(175.555,18.382)(180.050,15.979)
	(184.928,14.500)(190.000,14.000)
\drawline(190.000,14.000)(195.212,14.328)(200.254,15.690)
	(204.922,18.032)(209.029,21.259)(212.408,25.240)
	(214.924,29.817)(216.475,34.804)(217.000,40.000)
	(216.475,45.196)(214.924,50.183)(212.408,54.760)
	(209.029,58.741)(204.922,61.968)(200.254,64.310)
	(195.212,65.672)(190.000,66.000)
\drawline(187.000,10.000)(192.014,9.940)(196.969,10.713)
	(201.728,12.297)(206.157,14.648)(210.136,17.701)
	(213.553,21.371)(216.314,25.557)(218.343,30.143)
	(219.583,35.003)(220.000,40.000)(219.583,44.997)
	(218.343,49.857)(216.314,54.443)(213.553,58.629)
	(210.136,62.299)(206.157,65.352)(201.728,67.703)
	(196.969,69.287)(192.014,70.060)(187.000,70.000)
\put(215,53){\circle*{6}}
\put(216,27){\circle*{6}}
\dashline{4.000}(64,53)(107,53)
\dashline{4.000}(64,27)(107,27)
\dashline{4.000}(214,53)(257,53)
\dashline{4.000}(214,27)(257,27)
\drawline(0,80)(0,0)
\drawline(187,10)	(184.312,10.168)
	(181.998,10.288)
	(178.352,10.384)
	(174.000,10.000)
\drawline(174,10)	(170.803,9.608)
	(166.931,8.823)
	(163.093,7.876)
	(160.000,7.000)
\drawline(160,7)	(155.558,4.780)
	(151.934,2.805)
	(147.000,0.000)
\drawline(187,70)	(184.312,69.832)
	(181.998,69.712)
	(178.352,69.616)
	(174.000,70.000)
\drawline(174,70)	(170.803,70.392)
	(166.931,71.177)
	(163.093,72.124)
	(160.000,73.000)
\drawline(160,73)	(155.558,75.220)
	(151.934,77.195)
	(147.000,80.000)
\put(127,38){\makebox(0,0)[lb]{\smash{{{\SetFigFont{12}{14.4}{rm}=}}}}}
\put(84,40){\makebox(0,0)[lb]{\smash{{{\SetFigFont{20}{24.0}{rm}.}}}}}
\put(84,33){\makebox(0,0)[lb]{\smash{{{\SetFigFont{20}{24.0}{rm}.}}}}}
\put(234,40){\makebox(0,0)[lb]{\smash{{{\SetFigFont{20}{24.0}{rm}.}}}}}
\put(234,33){\makebox(0,0)[lb]{\smash{{{\SetFigFont{20}{24.0}{rm}.}}}}}
\put(234,47){\makebox(0,0)[lb]{\smash{{{\SetFigFont{20}{24.0}{rm}.}}}}}
\end{picture}
 }
\caption[]{}
\label{Kirby}
\end{figure}

\begin{figure}[htbp]
\begin{picture}(0,0)%
\epsfig{file=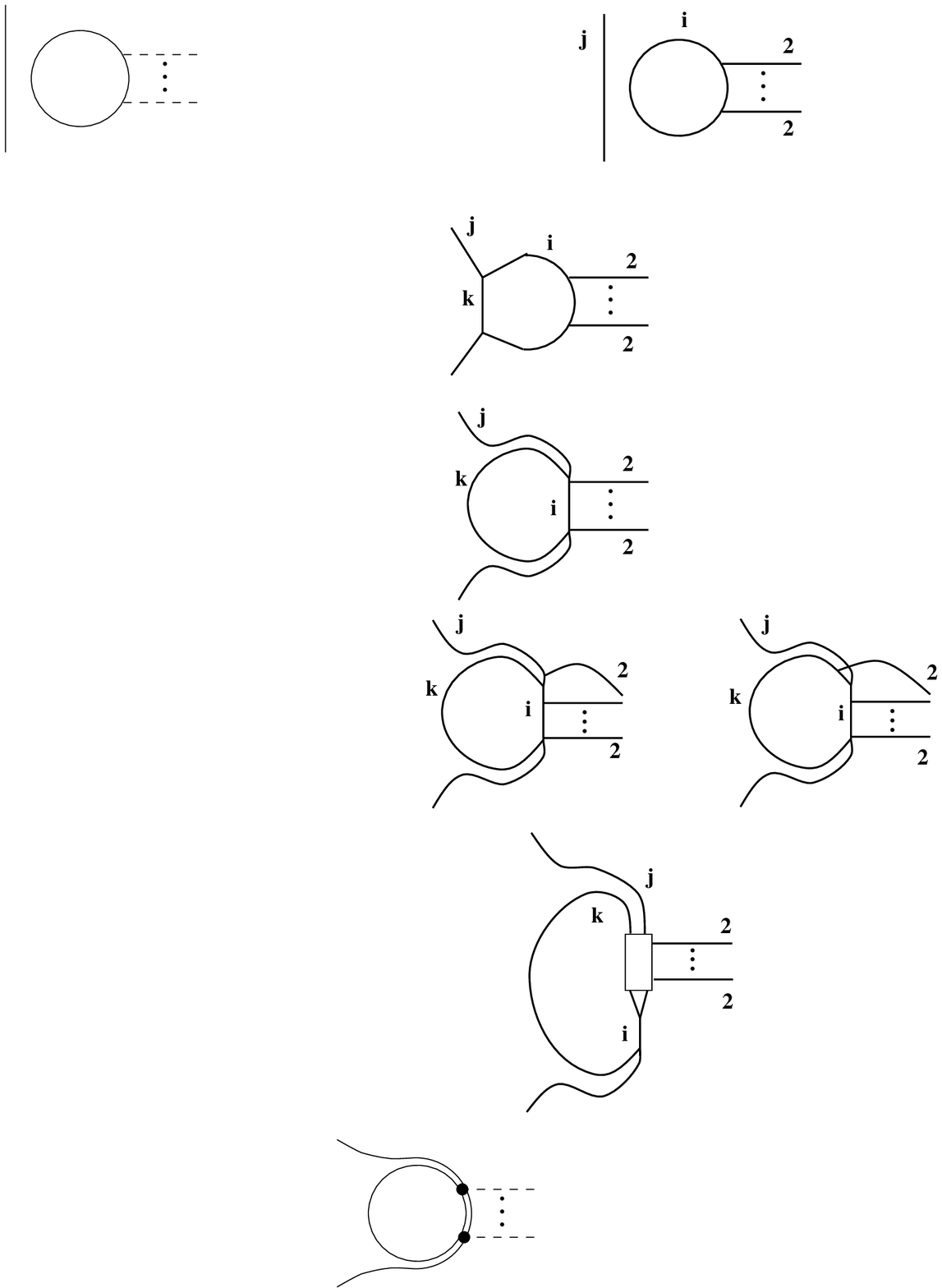}%
\end{picture}%
\setlength{\unitlength}{0.00083300in}%
\begingroup\makeatletter\ifx\SetFigFont\undefined
\def\x#1#2#3#4#5#6#7\relax{\def\x{#1#2#3#4#5#6}}%
\expandafter\x\fmtname xxxxxx\relax \def\y{splain}%
\ifx\x\y   
\gdef\SetFigFont#1#2#3{%
  \ifnum #1<17\tiny\else \ifnum #1<20\small\else
  \ifnum #1<24\normalsize\else \ifnum #1<29\large\else
  \ifnum #1<34\Large\else \ifnum #1<41\LARGE\else
     \huge\fi\fi\fi\fi\fi\fi
  \csname #3\endcsname}%
\else
\gdef\SetFigFont#1#2#3{\begingroup
  \count@#1\relax \ifnum 25<\count@\count@25\fi
  \def\x{\endgroup\@setsize\SetFigFont{#2pt}}%
  \expandafter\x
    \csname \romannumeral\the\count@ pt\expandafter
  \csname #3\endcsname}%
\fi
\fi\endgroup
\begin{picture}(8190,10460)(46,-10344)
\put(2701,-9811){\makebox(0,0)[lb]{\smash{\SetFigFont{14}{16.8}{rm}$ \omega^{(p)}  ( $}}}
\put(1576,-9811){\makebox(0,0)[lb]{\smash{\SetFigFont{14}{16.8}{sf}$=\cdots =$ }}}
\put(5041,-9811){\makebox(0,0)[lb]{\smash{\SetFigFont{12}{14.4}{sf}$)$}}}
\put(8236,-5786){\makebox(0,0)[lb]{\smash{\SetFigFont{12}{14.4}{sf}$)$}}}
\put(5746,-5756){\makebox(0,0)[lb]{\smash{\SetFigFont{12}{14.4}{rm}$ \displaystyle + i^{n-1} k $}}}
\put( 91,-5776){\makebox(0,0)[lb]{\smash{\SetFigFont{12}{14.4}{rm}$ \displaystyle = \sum_{<i,j,k>_p}^{p-2}  {(-1)^{i+j+k} (i+1) (j+1)(k+1) \over \Theta(i,j,k)}\, (i^{n-1}j$}}}
\put(781,-7868){\makebox(0,0)[lb]{\smash{\SetFigFont{12}{14.4}{rm}$ \displaystyle = \sum_{<i,j,k>_p}^{p-2}  {(-1)^{i+j+k} (i+1) (j+1)(k+1) \over \Theta(i,j,k)}$}}}
\put( 46,-616){\makebox(0,0)[lb]{\smash{\SetFigFont{14}{16.8}{rm}$ \omega^{(p)}  ( $}}}
\put(2311,-601){\makebox(0,0)[lb]{\smash{\SetFigFont{12}{14.4}{sf}$)$}}}
\put(2686,-601){\makebox(0,0)[lb]{\smash{\SetFigFont{12}{14.4}{rm}$ \displaystyle = \sum_{i,j=0}^{p-2}  {(-1)^{i+j} (i+1) (j+1)\, i^n}$}}}
\put(631,-2366){\makebox(0,0)[lb]{\smash{\SetFigFont{12}{14.4}{rm}$ \displaystyle = \sum_{<i,j,k>_p}  {(-1)^{i+j+k} (i+1) (j+1)(k+1) \over \Theta(i,j,k)}\,i^n$}}}
\put(631,-4031){\makebox(0,0)[lb]{\smash{\SetFigFont{12}{14.4}{rm}$ \displaystyle = \sum_{<i,j,k>_p}^{p-2}  {(-1)^{i+j+k} (i+1) (j+1)(k+1) \over \Theta(i,j,k)}\,i^n$}}}
\end{picture}

\caption{}
\label{proof}
\end{figure}

We are now able to state one of  the main theorems of this paper~:
Let $D$ be a chord diagram with support $X$, let
 $(\lambda_{i})i=1,\cdots, n$ be the coloring of the connected components 
 of $X$ and denote $\Delta_{\vec {\lambda}}= \prod_{i=1}^n 
 {(-1)^{\lambda_{i}}(\lambda_{i}+1)}$.
\begin{theo}
The weight system $\omega_{sl(2)}^{{(p)}}$ 
(resp. $\omega_{so(3)}^{(p)}$) given by~:
\begin{equation}
\omega^{{(p)}} (D)\=p \sum_{\vec{\lambda}} 
\Delta_{\vec{\lambda}} \omega^{T}_{\vec {\lambda}}(D),
\end{equation}
where the sum is over all
$\lambda_{i} \in I_{p} $ (resp.  $\lambda_{i} \in I_{p} \cap 2{\Bbb{Z}}$),
is a Kirby weight system valued into the field 
${{\Bbb{Z}} / p{\Bbb{Z}}}$.

\end{theo}
{\bf Proof} The proof is given in figure \ref{proof} where we use the 
formula \ref{identity} and the relation of figure \ref{weightint}.
The proof is given for the $sl(2)$ case, for the $so(3)$ case it is 
exactly the same proof taking into account 
 that if $(i,j,k)$ is an admissible 
triple and $i,j$ are even then $k$ is also even.
$\Box$

\begin{prop}
\label{integer}
Let $D\in {\cal A}((S^1)^L)$, let $\lambda_i$  be the highest
weight associated with the $i^{\mathrm{th}}$ component of $(S^1)^L$
and let $n_i(D)$ be the number of univalent vertices on the 
$i^{\mathrm{th}}$ component of $(S^1)^L$.
Then expanding $\omega^{{sl(2)}}_{\vec{\lambda}}$ in powers of $(\lambda_i+1)$ we get~:
\begin{equation}
\begin{array}{c}
{\displaystyle
\omega^{{sl(2)}}_{\vec{\lambda}}(D) =
\sum_{k_1 \geq 0}^{[{n_1\over 2}]}\cdots \sum_{k_L \geq 0}^{[{n_L\over 2}]}
{2^{-d^{\circ}(D)} \over \prod_{i=1}^L (n_i(D)+1)!} a_{k_1, \cdots, k_L}(D) 
(\lambda_1+1)^{2k_1+1}\cdots (\lambda_n+1)^{2k_L+1} }\\
{\displaystyle
=\sum_{k_i \geq 0}^{[{n_i(D)\over 2}]}(\lambda_1+1)^{2k_1+1}
2^{-d^{\circ}(D)} a_{\vec{\lambda},K_i}(D) }
\end{array}{c}
\end{equation}
where $d^{\circ}(D)$ is the degree of the diagram $D$ and $[{n_1\over 2}]$ denotes the integral part of ${n_1\over 2}$.
Moreover, the coefficients $a_{k_1, \cdots, k_L}(D)$ lies in $\Bbb{Z}$, and 
$a_{\vec{\lambda},k_i}(D)$ lies in $ {\Bbb{Z}}$ if $\lambda_j \in {\Bbb{Z}}$.
\end{prop}

{\bf Proof}
Let ${\cal U}(sl(2))_{\Bbb{Z}}$ the subring of 
${\cal U}(sl(2))$ generated over $\Bbb{Z}$ by 
the Chevalley basis elements $e,f,h$ (\ref{sl2}).
The structure constants appearing in the commutation relations of the 
Chevalley generators are integers; thus 
${\cal U}(sl(2))_{\Bbb{Z}}$ admits a $\Bbb{Z}$-basis consisting of 
all $ e^ah^bf^c, a,b,c \in \Bbb{N}$.

Using the definition of the $sl(2)$ weight system and the fact 
that $2C \in \cal{G}_{\Bbb{Z}} \otimes \cal{G}_{\Bbb{Z}}$ where 
${\cal G}_{\Bbb{Z}}$ 
denotes the $\Bbb{Z}$-linear span of the Chevalley basis, it is clear 
that $2^{d^{\circ}(D)}\omega^{sl(2)}_{\vec{\lambda}}(D)$ is 
constructed as the trace over $V_{\lambda_1} \otimes \cdots \otimes 
V_{\lambda_l}$, $\lambda_i \in \Bbb{N}$ of an element $ x(D) \in 
({\cal U}(sl(2))_{\Bbb{Z}})^{\otimes L})^{\cal G}$, where the power of 
$\cal G$ means that $x(D)$ commutes with the diagonal action of $\cal 
G$.  Moreover the degree of $x(D)$ with respect to the 
i$^{th}$-component of $({\cal U}(sl(2))_{\Bbb{Z}})^{\otimes L})$ is at 
most $n_i(D)$.

Let $ Z(sl(2))$  be the center of ${\cal U}(sl(2))$
and let $V_\lambda$ the irreducible highest weight module of 
(non necessarily integral) weight $\lambda$. 
If $z \in  Z(sl(2))$, $z$ acts as a scalar on $V_\lambda$, 
denoted as $\chi_\lambda(z)$.  From Harish-Chandra's theorem we know that 
the infinitesimal character $\chi_\lambda(z)$ is an even polynomial in 
$\lambda +1$ whose degree is at most the degree of $z$ with respect to 
the filtration of ${\cal U}(sl(2))$.  In the case of $sl(2)$ the 
center $Z(sl(2))$ is generated by $C$ and 
$\chi_\lambda(C)={((\lambda+1)^2-1)\over2}$.  Let ${\cal{D}}=[{\cal 
U}(sl(2)),{\cal U}(sl(2))]$ be the subspace of ${\cal U}(sl(2))$ 
generated by all commutators.  We then have that ${\cal 
U}(sl(2))=Z(sl(2))\oplus \cal{D}$ (\cite{Hum}) so we can extend 
$\chi_\lambda $ to all of ${\cal U}(sl(2))$ by requiring it to be $0$ 
on $\cal{D}$.  When $\lambda \in \Bbb{N}$, $V_\lambda$ is finite 
dimensional and $tr_{V_\lambda}(x)=(\lambda+1)\chi_\lambda(x)$.  Let 
$v_i, \, i=0,\cdots,\lambda$ the basis of $V_\lambda$ given in 
(\ref{sl2}) and denote ${V_\lambda}_{\Bbb{Z}}$ the $\Bbb{Z}$-span of 
$v_i$, if $x \in {\cal U}(sl(2))_{\Bbb{Z}}$ then $x\cdot v_i \in 
{V_\lambda}_{\Bbb{Z}}$ so $tr_{V_\lambda}(x) \in {\Bbb{Z}}$.

We can conclude from this analysis that~:
$2^{d^{\circ}(D)}\omega^{sl(2)}_{\vec{\lambda}}(D)=
{\displaystyle \prod_{i=1}^{L}}(\lambda_i +1) 
\chi_{\lambda_1}\otimes \cdots\otimes\chi_{\lambda_L} x(D)$
is an odd polynomial of degree $n_i(D) +1$ in $\lambda_i+1$ taking
integer values when  $\lambda_i,i=1,\cdots,L$ are integers.
It is a standard exercise to show by recurrence 
that such a polynomial is a 
${\Bbb{Z}}$-linear combination of the polynomials
$$
\displaystyle{\prod_{i=1}^L }
\left( { \displaystyle{\lambda_i +1+k} 
\atop {\displaystyle{2k_i+1}}}\right),\, 2k_i\leq n_i(D),
$$
where ${\displaystyle \left( { x 
\atop k}\right)= {x(x-1)\cdots(x-k+1) \over k!}}$. 
Thus $2^{d^{\circ}(D)}{\prod_{i=1}^L (n_i(D)+1)!} 
\omega^{sl(2)}_{\vec{\lambda}}(D)$
is a polynomial with integer coefficients.

In order to prove the second part of the theorem
we consider the partial trace 
$ x_1(D)={\displaystyle \prod_{i=2}^{L}}(\lambda_i +1) 
1\otimes \chi_{\lambda_2}\otimes \cdots\otimes\chi_{\lambda_L} x(D)$ with 
 $\lambda_2,\cdots \lambda_L \in \Bbb{N}$.
This is a well defined operation ( i-e which do not depend on the particular choice
 of $x(D)$ ) moreover
$x_1(D)\in {\cal U}(sl(2))_{\Bbb{Z}} \cap Z(sl(2)) $ thus 
$x_1(D) $ is a polynomial of degree $n_i(D)/2$ in $2C$ with integral
coefficients. \\
$\Box$

{\bf Preliminary}

Consider $D_{< 2n}$ the subspace of ${\cal{A}}(X)$ 
generated by all chord diagrams which possess less than $2n$ univalent 
vertices on one component of $X$, and denote$ by \phi_n: {\cal{A}} 
(X)\rightarrow {\cal{A}} (X)/D_{< 2n} $ the canonical projection.  By 
definition,
\begin{eqnarray}
\phi_n (D) =0 \,{\mathrm{ if }} \, n_i(D) < 2 n,
\end{eqnarray}
where $n_i(D)$ denote the number of univalent vertices on the 
$i^{\mathrm{th}}$ component of the diagram $D$.
We define $N_p = ({p-3 \over 2})$ and we denote by 
${\Bbb{Z}}_{(n)}\equiv {\Bbb{Z}}[{1\over 2}, \cdots, {1\over n}]$ the 
ring generated over ${\Bbb Z}$ by $\{{1\over 2}, \cdots, {1\over 
n}\}$.
\begin{prop}

Let $D \in {\cal{A}} ((S^1)^L)\otimes {\Bbb{Z}}_{(p-1)}$, then~:
\begin{equation}
\omega_{sl(2)}^{(p)}(D) {\=p}
 2^{L}\, \omega_{so(3)}^{(p)}(D)
,\, if\,  n_i(D)\leq 2N_p +1
\end{equation}
\centerline{$ \omega_*^{(p)}(D) 
\=p \omega_*^{(p)} \circ \phi_{N_p}(D).$}
\label{propo} 
\end{prop}
where $*$ in $\omega_*$ means $sl(2)$ or $so(3)$.
This  proposition  state that $sl(2)$ and $so(3)$ Kirby weight
 systems are essentially the same.

Define 
\begin{equation}
\epsilon_p(i) \=p \left\{ 
\begin{array}{cc}
-1, & \mathrm{if}\, p-1\,\mathrm{divides} \, i, \\
0, &\,\mathrm{if\, not}
\end{array}
\right. .
\end{equation}
We have the following lemma~:
\begin{lemma}
\label{Vst}
\begin{equation}
\begin{array}{c}
{\displaystyle
{\sum_{k=0}^{p-2}}  (k+1)^i  \=p \epsilon_p(i) }\\
{\displaystyle
{\sum_{k=0\atop k\equiv 0 mod(2)}^{p-2}} (k+1)^{2i}
\=p {1\over 2} 
\epsilon_p(2i) }
\end{array}
\end{equation}
\end{lemma}
 
The first equality, which is the Von-Staudt theorem, follows directly  
from the fact that the group of invertible elements of ${\Bbb{Z}}/p{\Bbb{Z}}$ 
 is the cyclic group of order $p-1$ for $p$ prime.  
 for the second 
 equality consider~: 
 $${ \epsilon_p(2i)\=p {\sum_{k=1}^{p-1}} k^{2i} \=p 
 {\sum_{k=1}^{p-1\over 2}}k^{2i} + \sum_{k={p-1\over 2}+1}^{p-1} 
 k^{2i}},$$ making the change of variables $j =p - k$ the second term 
 in the LHS can be written as ${ {\sum_{j=1}^{p-1\over 
 2}}(p-k)^{2i}\=p {\sum_{j=1}^{p-1\over 2}}(k)^{2i}.}$ Thus ${ 
  {\sum_{k=1}^{p-1\over 2}}k^{2i}\=p{1\over2} \epsilon_p(2i).}$ 
 Together with, $${ \epsilon_p(2i)\=p {\sum_{ k=1}^{p-1\over 2}} 
 (2k-1)^{2i} + {\sum_{ k=1}^{p-1\over 2}} (2k)^{2i}, }$$ this imply 
 $${\displaystyle {\sum_{ k=1}^{p-1\over 2}} (2k-1)^{2i} 
 \=p(1-2^{2i-1}) \epsilon_p(2i)
}.$$
We get the desired result using the fact that $ 2^{2i} \=p 1$ 
if $p-1$ divides $2i$.

{\bf Proof of Proposition \ref{propo}}\\
Let $\psi_p:{\Bbb{Z}}_{(p-1)} \rightarrow {{\Bbb{Z}}\over p{\Bbb{Z}}}$ 
be the evaluation map.
By definition we have for $D \in {\cal{A}}((S^1)^L)$. 
\begin{equation}
\omega^{(p)}_{sl(2)}(D)
\=p \psi_p \left( \sum_{\lambda_1=0}^{p-2} \cdots \sum_{\lambda_L=0}^{p-2} 
\prod_{i=1}^L (\lambda_i +1) \omega_{\vec{\lambda}}^{sl(2)}(D) \right)
\end{equation}
This is well defined since we have seen in prop.\ref{integer} that 
$\omega_{\vec{\lambda}}^{sl(2)}(D) \in {\Bbb{Z}}_{(2)}$.
If $D$ is such that $ n_i(D)\leq 2N_p +1$ this means that 
the coefficients of $ \prod_{i=1}^{L}(\lambda_i +1)^{2(k_i+1)}$ in the development
 of $\omega_{\vec{\lambda}}(D)$ belongs to 
${\Bbb{Z}}_{(p-1)}$, thus we can distribute $\psi_p$ and apply 
the preceding 
lemma to get the first conclusion of the proposition.

If $D$ is such  that $\phi_{N_p}(D)=0$, for example $n_1(D)<2 N_p$, 
then using the results of prop.\ref{integer} especially the fact that 
the coefficient of $(\lambda_1 +1)^{2(k_1+1)}$ belong to 
${\Bbb{Z}}_{(2)}$ we have 
\begin{equation}
\omega^{(p)}_{sl(2)}(D) \=p 
\sum_{k_1 \geq 0}^{[{n_1\over 2}]} \sum_{\lambda_1=0}^{p-2}
\psi_p({ (\lambda_1 +1)^{2(k_1+1)} }) \psi_p( a_{k_1}(D))
 \=p
\sum_{k_1 \geq 0}^{[{n_1\over 2}]}
\epsilon_p(2(k_1+1))\psi_p( a_{k_1}(D))
\end{equation}
 $k_1 < N_p$ imply $2(k_1+1)<p-1$ thus $\epsilon_p(2(k_1+1))\=p 0$.\\
$\Box$

\section{The asymptotic Invariant for rational homology 3-sphere}

Let  $L$ be a framed  link and ${\hat Z}(L)$ be the canonical 
Vassiliev invariant \cite{AF1,LeM1,Car}. Which is a formal power series 
in $\hbar$~: ${\hat Z}(L)={\sum_{m=0}^{+\infty}\hbar^m{\hat Z}_m(L)}$ 
with ${\hat Z}_m(L)\in {\cal A}_m(L) \otimes {\Bbb{Q}}$.  Let $\nu$ be 
the value of $\hat Z$ for the un-framed trivial knot and define 
${\check Z}(L)=\nu \otimes \cdots \otimes \nu \cdot {\hat Z}(L) $, 
where the product is realised with the connected sum along each 
component of $L$.  One of the main theorem of \cite{LeM2} states that 
if $\omega $ is a Kirby weight system then $\omega\circ {\check Z}(L)$ 
is invariant under the second Kirby move (hand-slide) and under the 
change of orientation.  Once we have a Kirby weight system we can 
construct, up to a normalisation problem, a framed link invariant, 
independent of the orientation and invariant under all Kirby moves, 
hence a 3-manifold invariant.  The problem we face in our case comes 
from the fact that the weight systems we constructed are valued into a 
non-zero characteristic field.  In order to apply our weight system to 
the universal invariant we need to know the integrality properties of 
its coefficients.  T.Q.~Le in a recent paper \cite{Le3} completely 
study the denominators of the Kontsevich integrals.  Combining the 
propositions 5.3 and 6.1 of \cite{Le3}, we can state~:
 \begin{lemma}
\label{integr}
If $L$ is an algebraically split (ASL) link then~:
\begin{equation}
\phi_n({\check Z}_m(L)) \in {\cal A}_m(L)\otimes {\Bbb{Z}}_{(2(n+1))},
\,{ \mathrm{for}} \,  m\leq n(|L|+1)
\end{equation}
where $|L|$  denotes the number of connected
component of the link $L$.
\end{lemma}
$L$ is called algebraically split if it is a framed link with diagonal 
linking matrix.
It is known that if $M$ is rational homology 3-sphere then there exist 
Lens spaces such that the connected sum of $M$ with these lens spaces 
can be obtain by surgery along an algebraically split link possessing 
a non degenerate linking matrix \cite{LeM2}.  Using the multiplicative 
property of the invariant under connected sum this means that the 
computation of the asymptotic invariant we are going to make 
is applicable to atleast the case of all rational homology 3-spheres.   
\begin{defi}
If $L$ is an ASL we can define~: 
\begin{equation}
F_*^p(L)= \sum_{m=0}^{N_p}\hbar^m \omega^{(p)}_*({\check Z}_{N_pL+m})(L)
= \hbar^{-N_p|L|}\omega^{(p)}({\check Z}(L) \,\,
{\mathrm{ mod}}(\hbar^{N_p(|L|+1)+1})),
\end{equation}
where $p$ is a prime odd integer and $N_{p}=(p-3)/2$.
\end{defi}
$\omega^{(p)}$ here means $\omega^{(p)}\circ \psi_{N_p}$,
which are by proposition \ref{propo}  identical on ${\cal A}(L)
\otimes {\Bbb{Z}}_{(p-1)}$, and $*$ refers to $sl(2)$ or $so(3)$.
 $\omega^{(p)}$ is well defined on 
${\cal A}(L) \otimes {\Bbb{Z}}_{(p-1)}$, thus from lemma \ref{integr} 
the definition of $F^{p}$ suffers from no ambiguity.

Let $U^{\pm} $ be the unknot with framing $\pm 1$, $F_p(U^{\pm})$ 
is an invertible element (i-e the coefficient of $\hbar^0$ 
is non zero, see \ref{U+}).  For $L$ an ASL we can define~:
\begin{equation}
O_*^p(L) = {F_*^p(L) \over F_*^p(U^{+})^{\sigma_+} F_*^p(U^{-})^{\sigma_-} }.
\end{equation}
where $\sigma_+$ (resp.  $\sigma_-$)  is the number of positives
 (resp. negatives) eigenvalues of the linking matrix of $L$.
 $O_*^p(L) $ is invariant under all Kirby moves.

\subsection{Fermat limit}

For a complete discussion of the notion of Fermat limit, see 
\cite{Lin}.  For the reader's convenience we recall here some 
definitions that we are going to use.

\begin{defi}
Let  $ (u_p)_{p\, \mathrm{prime}} $ be a sequence of 
${{\Bbb{Z}} /p{\Bbb{Z}}}$ numbers.
We say that $ (u_p)_{p\, \mathrm{prime}} $ admits a Fermat-limit if 
there exists a rational number $u$ and $N \in \Bbb{N}$ such that 
$\psi_p(u) \=p u_p $ for all $p > N$.
 Consider 
$$ \displaystyle{ {\cal{S}}_p=\{ F_p =\sum_{n=0}^{N_p} \hbar^n F_{n,p},
F_{n,p} 
\in {{\Bbb{Z}} \over p{\Bbb{Z}}}\}}.$$
We say that $(F_p)_{p \, \mathrm{prime}}$, $ F_p \in{\cal{S}}_p$ admits 
a strong Fermat limit if there exist $F \in {\Bbb{Q}}[[\hbar]]$, 
$\displaystyle{ F =\sum_{n=0}^{+\infty} \hbar^n F_{n} }$ and $N \in 
\Bbb{N}$, such that $\displaystyle{ \psi_p (F \, 
{\mathrm{mod}}(\hbar^{N_p+1}) ) \equiv \sum_{n=0}^{N_p} \hbar^n 
\psi_p(F_{n}) \=p F_p }$ for all $p > N$.
\end{defi}
\begin{defi}
Given integers $f_1,\cdots,f_n$ we consider a formal integration 
procedure defined by~:
\begin{equation}
\label{int}
I_{\vec{f},\hbar} ( \prod_{i=1}^{|L|} (\hbar \alpha_i)^{k_i}) =
\left\{
\begin{array}{cc}
\displaystyle{\prod_{i=1}^{|L|} \left({-\hbar \over f_i}\right)^{k_i  \over 2} 
{k_i ! \over({k_i  \over 2})!}},
  &{\mathrm{if}}\, k_i {\mathrm{\, is \, even }}
 \forall i \in \{1,\cdots ,|L|\}\\
0,& \,\mathrm{if\, not}  
\end{array}
\right. .
\end{equation}
We can  extend by linearity $I_{\vec{f},\hbar}$ to all formal power
 series of the form
$\sum_{n,{\vec k}} S_{n,{\vec k}} 
\hbar^n (\prod_{i=1}^{|L|} \alpha_i^{2k_i})$
if $\{S_{n,{\vec k}}\neq 0| \, n+\sum_i k_i=q\} $ is a finite set
for all $q$.
\end{defi}
{\bf Remark} \\
If we consider $ \hbar $ as an imaginary number whose value is
${2i \pi \over K}$ then $I_{\vec{f},{2i \pi \over K}}$ is really an integral~:
\begin{equation}
I_{f,{2 i\pi \over K}}(P(\alpha))=
e^{-{i\pi \over 4} sgn(f)} \left( {2K \over |f|}\right)^{1\over 2}
lim_{\epsilon\rightarrow 0} \int_{-\infty}^{+\infty}d\alpha\, 
e^{{i\pi\over 2K} f\alpha^2} P(\alpha)e^{-\epsilon\alpha^2}
\end{equation}
\begin{theo}
\label{strongFer}
Consider $N,n\in \Bbb{N}$ and 
$D \in {\cal A}_n(L)\otimes {\Bbb{Z}}_{(N)}$
Then, 
\begin{equation}
\left({ (N_p+1)! \over \epsilon(*) \hbar^{N_p}}\right)^{|L|}
 {1\over \prod_{i=1}^{|L|}\left( {f_i\over p} \right) }
 \omega^{(p)}_*(
\prod_{i=1}^{|L|}e^{ {\hbar f_i \over2}(\theta_i )}
 \cdot \hbar^{d^{\circ}(D)}D \, \mathrm{mod}(\hbar^{N_p(L+1)+1})
), 
\in {\cal{S}}_p
\end{equation}
admits a strong Fermat limit which is 
\begin{equation}
(\prod_{i=1}^{|L|}e^{\hbar f_i \over 4 })
\hbar^{|L|}I_{\vec{f},\hbar}\left(
\prod_{i=1}^{|L|} \alpha_i \tilde{\omega}_{\vec{\alpha}}
(\hbar^{d^{\circ}(D)}D)
\right),
\end{equation}
where $*$ denotes $sl(2)$ or $so(3)$,
$\epsilon(sl(2))=-1, \epsilon(so(3)) =-1/2$,
$\left( {\cdot \over p} \right)$ is the Legendre symbol and 
$\tilde{\omega}_{\vec{\alpha}} = \omega_{(\alpha_1-1,\cdots, 
\alpha_{|L|}-1)}$.

Moreover, if $D$ is a chord diagram,
$\hbar^{|L|}I_{\vec{f},\hbar}\left(
\prod_{i=1}^{|L|} \alpha_i \tilde{\omega}_{\vec{\alpha}}(
\hbar^{d^{\circ}(D)} D)
\right)$ is a polynomial in $\hbar$ which degree is lower than 
$ d^{\circ}{(D)}$ and whose valuation is greater than half the number of 
internal vertices of $D$.
\end{theo}
Note that the Fermat limit does not depend on the choice of $sl(2)$ 
or $so(3)$.

\begin{theo}
\label{Fer}
Letting $F_{n}^p$ be the coefficient of $\hbar^n$ in the developpement of 
$F^p$, then if $L$ is an ASL such that the framings $f_i$ of all its 
components are non zero then
\begin{equation}
{  ((N_p+1)!)^{|L|} \over \prod_{i=1}^{|L|}\left({f_i \over p}\right) }
\, F_{n}^p(L)
\end{equation}
admits a Fermat limit denoted $F_n$.
Moreover, $F=\sum_{n=0}^{+\infty}\hbar^n F_n$ is equal to~:
\begin{equation}
F =\prod_{i=1}^{|L|}e^{-{\hbar f_i \over 4} }
\hbar^{|L|}I_{\vec{f},\hbar}\left(
\prod_{i=1}^{|L|} \alpha_i
 \tilde\omega_{\alpha}({\check Z}(L)) \right)
\end{equation}
\end{theo}

For example 
\begin{equation}
\label{U+}
F(U^{\pm}) = 
{\mp 2 \hbar \over e^{ \hbar \over 2} - e^{-{ \hbar \over 2}} }
 e^{\mp{ 3\hbar \over 4}}
\end{equation}

The asymptotic rational homology 3-sphere quantum invariant $O$
 is defined as the Fermat limit of $O_p$.
If we denote $J_{\vec{\alpha}}$ the colored Jones polynomial then~:
\begin{equation}
O(L)= e^{-{\hbar\over 4}\sum_i(f_i +{1\over f_i} -{3\over4} sgn(f_i))}
 \prod_{i=1}^{|L|}sgn(f_i)
I_{\vec{f}, \hbar}(J_{\vec{\alpha}+{1\over \vec{f}}}(L))
\end{equation}
This expression correspond to the Rozansky formulation of the 
 Ohtsuki invariant 
\cite{Roz}.

\subsection{Proofs}

Let $ D \in {\cal A}_m(L)\otimes {\Bbb{Z}}_{(N)}$, $N\in \Bbb{N}$ and denote
by $n_i(D)$ the number of chords arriving the i$^{th}$ component of 
$L$.  Let $\theta_i \in {\cal A}(L) $ be the diagram possessing only 
one chord on the i$^{th}$ component of $L$.

Using the definitions of mod$(\hbar^{N_p(L+1)+1})$ and  $\phi_{N_p}$
we have the expansion,
\begin{equation}
\phi_{N_p}(\prod_{i=1}^{|L|}e^{{\hbar f_i \over2}(\theta_i ) }
\cdot \hbar^{d^{\circ}(D)} D \, \mathrm{mod}(\hbar^{N_p(L+1)+1}))
= \sum_{\vec l \in {\cal{J}}_{p}(D)}
\prod_{i=1}^{|L|}\left({\hbar f_i \over2}\right)^{l_i}
{(\theta_i)^{l_i} \over l_i!}
\cdot \hbar^{d^{\circ}(D)} D,
\end{equation}
 where ${\cal{J}}_{p}(D)=\{\vec l \in
 {\Bbb{N}}^{|L|}, d^{\circ}(D) + \sum_i l_i \leq N_p(|L|+1),\, n_i(D)
 + 2l_i \geq 2N_p \}$.
Remembering that $2d^{\circ}(D)$ is the number of vertices of the 
diagram $D$ we have $ 2N_p \geq 2(d^{\circ}(D)-N_p|L|) + \sum_i 2l_i 
\geq \sum_i (n_i(D)+2l_i -2N_p)$, which implies that $n_i(D)+2l_i\leq 4N_p$, 
and a fortiori $l_i \leq 2N_p$.

Thus $\left({\theta_i ^{l_i} \over l_i!}
\cdot D\right)\in{\cal A}(L)\otimes {\Bbb{Z}}_{(max(N,2N_p))} $ if 
$\vec l \in {\cal{J}}_{p}(D)$, so we can apply $\omega^{(p)}$ to this 
diagram for $p>N$. 
Using the definition of $\omega^{(p)}$, $p>N$ we have to compute~:
\begin{equation}
\sum_{\vec{\alpha}=\vec{1}}^{\vec{p}-\vec{1}}
\sum_{\vec l \in  {\cal{J}}_{p}(D) }
 \psi_p \left( 
\prod_{i=1}^{|L|}\left({\hbar f_i \over2}\right)^{l_i}\alpha_i 
\,\tilde{\omega}_{\vec \alpha} \left( 
\prod_{i=1}^{|L|}{\theta_i^{l_i} \over l_i!}\cdot \hbar^{d^{\circ}(D)} D
\right)\right).
\end{equation}

Using the notation of proposition \ref{integer} and the fact that $$
\tilde{\omega}_{\vec \alpha} \left( \prod_{i=1}^{|L|} {\theta_i^{l_i} 
\over l_i!}\cdot D \right) = \prod_{i=1}^{|L|} ({\alpha_i^2-1 \over 
2})^{l_i} {1\over l_i!}
\tilde{\omega}_{\vec \alpha} \left(\hbar^{d^{\circ}(D)}D \right) $$
 this can be expressed as~:
\begin{equation}
\label{dev}
\sum_{\vec{k}\geq 0 \atop 2\vec{k} \leq \vec{n}(D) }
\sum_{\vec{\alpha}=\vec{1}}^{\vec{p}-\vec{1}}
\sum_{\vec l \in  {\cal{J}}_{p}(D) } 
 \psi_p \left( \hbar^{d^{\circ}(D)}a_{\vec{k}}(D) 
\prod_{i=1}^{|L|}\left({\hbar f_i \over2}\right)^{l_i}
{\alpha_i}^{2(k_i+1)}({\alpha_i^2-1 \over 2})^{l_i} {1\over l_i!}
  \right).
\end{equation}
We know from proposition {\ref{integer} that 
$a_{\vec{k}}(D) \in {\Bbb{Z}}_{(max(n_i(D)+1,N)}$, thus if 
 $ p > max(n_i(D)+1,N) $  we can distribute $\psi_p$ over all 
the factors appearing in \ref{dev}.
Using the result of Lemma \ref{Vst} we get~:
\begin{equation}
\psi_p \left(
\prod_{i=1}^{|L|}
\sum_{{\alpha_i}=1}^{{p}-{1}}
{\alpha_i}^{2(k_i+1)}({\alpha_i^2-1 \over 2})^{l_i} {1\over l_i!}
\right)
\=p 
\sum_{\vec{q}+\vec{r}=\vec{l}}
\prod_{i=1}^{|L|} {(-1)^{r_i} \over 2^{l_i} q_i ! r_i!}
\epsilon_p(2(k_i +q_i+1)).
\end{equation}
But $2(k_i +q_i+1) \leq n_i(D) + 2 l_i +2 \leq 4N_p+2=2(p-1)-2$, 
so $\epsilon_p(2(k_i +q_i+1))\neq 0 $ if and only if $ k_i +q_i =N_p$.
The computation for the case of $so(3)$ is the same except that we have 
to replace $\epsilon_p$ by ${1\over 2} \epsilon_p$ 
(see proposition \ref{propo}).

Making the change of variable $ j_i = l_i +k_i -N_p $, 
we can express  (\ref{dev}) in the following way~: 
\begin{equation}
\sum_{\vec{k},\vec{j}}
\psi_p(a_{\vec{k}}(\hbar^{d^{\circ}(D)}D) )
\psi_p\left(
\prod_{i=1}^{|L|}\left({\hbar f_i \over 4}\right)^{j_i + N_p -k_i}
{(-1)^{j_i} \over j_i! (N_p -k_i)!}
\right),
\end{equation}
where the summation is over all $\vec{k},\vec{j}$ satisfying 
$ 0\leq 2k_i \leq  n_i(D), j_i \geq 0$ and 
$ d^{\circ}(D) + \sum_i (j_i -k_i) \leq N_p$
 (recall that we have suppose that $ p > max(n_i(D)+1,N)) $.

Using the identity~:
\begin{equation}
{1 \over (N_p -k)!} \=p {1\over(N_p+1)!} \left({-1 \over 4}\right)^{k+1} 
{(2(k+1))! \over (k+1)!},
\end{equation}
we get if $\psi_p(f_i) \neq 0 $
\begin{eqnarray}
\displaystyle{
\left({\epsilon(*)\hbar^{N_p} \over  (N_p+1)!}\right)^L 
\prod_{i=1}^{|L|}\left({f_i \over p}\right) \times } \\
\displaystyle{
 \left(\sum_{\vec{k}\geq 0 \atop 2\vec{k}\leq \vec{n}(D)}
\psi_p(a_{\vec{k}}(\hbar^{d^{\circ}(D)}D ) ) 
\hbar^{ |L|}\psi_p( 
I_{\vec{f}}(\prod_{i=1}^{|L|} \alpha_i^{2(k_i+1)}))
\right)
\psi_p(\prod_{i=1}^{|L|}e^{-\hbar f_i \over 4} )
\, mod(\hbar^{N_p(|L|+1)+1}),}
\end{eqnarray}
where we use the fact $4^{N_p+1}\=p 1$, 
$ f_i^{(N_p+1)}\=p\left({f_i \over p}\right)$.
This proves the first part of the theorem \ref{strongFer}.
We know from \ref{integer} that 
$\prod_{i=1}^{|L|} \alpha_i \tilde{\omega}_{\vec{\alpha}}(
\hbar^{d^{\circ}(D)} D)$ is a polynomial of total degree, with respect 
to $\alpha$,
smaller that $ \sum_{i=1}^{L}(n_i(D)+2)$, this means that the valuation 
of $ \hbar^{|L|}I_{\vec{f},\hbar}\left(
\prod_{i=1}^{|L|} \alpha_i \tilde{\omega}_{\vec{\alpha}}(
\hbar^{d^{\circ}(D)} D)
\right)$ is greater than 
${d^{\circ}(D)}+ L - (1/2)\sum_{i=1}^{L}(n_i(D)+2)=I(D)/2$,
where $I(D)$ denote the number of internal vertices of the graph 
$D$.\\
$\Box$

Let $D$ be a chord diagram, denote $\bar D$ the corresponding graph 
( the graph obtained by removing all Wilson circles from $D$).  We say 
that $D\in I_1$ if each connected component of the graph $ \bar D$ 
contains at least one trivalent vertices.  And an element $D\in 
{\cal{A}}_n \otimes {\Bbb{Q}}$ is said to be in $I_1$ if it can be 
expressed as a sum of diagrams belonging to $I_1$.
\begin{lemma}
If $D\in I_1 $ then $I(D) \geq {d^{\circ}(D) \over 2}$,
\end{lemma}
{\bf Proof}\\
If $D \in I_1$ is a connected graph then $d^{\circ}(D)\geq 2$ thus, if
 $c({\bar D})$ denotes the number of connected components of $\bar D$,
 then 
 $ d^{\circ}(D) \geq 2c({\bar D})$. 
Using the Euler characteristic of $\bar D$ and the fact  that the 
graph is trivalent we conclude that
$ I(D) +2c({\bar D}) \geq N(D)$, where $N(D)=\sum_{i=1}^{L}n_i(D)$
is the number of univalent  vertices of ${\bar D}$,  thus 
$2d^{\circ}(D)=N(D) +I(D)\leq 2I(D) +2c({\bar D})\leq  2I(D) + 
d^{\circ}(D)$.\\
$\Box$

If $L$ is an ASL, and $L^\prime $ is the asociated link of framing 
zero, then $\check Z(L)=\prod_{i=1}^{|L|}e^{{\hbar f_i \over2}(\theta_i )}
\check Z_n(L^\prime)$, where $\theta_i$ denotes the chord diagram 
possessing only one chord one the $i^{th}$-component and the product 
is realised as usual by the connected sum, moreover $\check 
Z_n(L^\prime) \in I_1$ \cite{Le3}.  This implies, by the second part of 
theorem \ref{strongFer} and the preceding lemma, that $ 
\hbar^{-N_p|L|}\omega^{(p)}(\hbar^n {\check Z}_n(L^\prime))$ has a 
valuation with respect to $\hbar$ greater than ${n \over 4}$.  We can 
conclude that if $n\leq N_p(L+1)$ then $F_{p,n}$ is equal to the 
coefficient of $\hbar^n$ in the development of~:
\begin{equation}
{(N_p+1)! ^{|L|} \over \prod_{i=1}^{|L|}\left(f_i \over p\right)}
\sum_{k=0}^{4n} \omega^{(p)}( 
\prod_{i=1}^{|L|}e^{{\hbar f_i \over2}(\theta_i )}
\check{Z}_{k}({L}^\prime)  ).
\end{equation}
From theorem \ref{Fer} each term of this finite sum admits a Fermat limit; 
hence $F_{p,n}$ admits a Fermat limit which is 
 the coefficient of 
$\hbar^n$ in the development of~:
\begin{equation}
F(L)=\prod_{i=1}^{|L|}e^{-\hbar f_i \over 4}
\hbar^{|L|}I_{\vec{f}}(\prod_{i=1}^{|L|}\alpha_i
\tilde{\omega}_{\vec{\alpha}}(\check Z (L)) ). 
\end{equation}
\\
$\Box$

Let us denote $U_\pm$ the unknot with framing $\pm$, $\nu $ the value 
of $ \hat {Z}$ on the un-framed unknot and $ J_{\vec{\alpha}}(L) 
=\omega_{\vec{\alpha}}(\hat{Z}(L^{\prime}))$ the colored Jones 
polynomial.  $$\omega_{{\alpha}}(\nu) = { e^{\hbar \alpha \over 2} - 
e^{-{\hbar \alpha \over 2}} \over e^{\hbar \over 2} - e^{-{\hbar \over 
2}}}$$ and
\begin{equation}
I_{f,\hbar} (e^{q\hbar\alpha} P(\alpha)) = (e^{- {q^2 \hbar\over f}}I_{f,\hbar}( P(\alpha- {2q\over f})).
\end{equation}

Thus 
\begin{equation}
F(L) = { -2 \hbar \over e^{\hbar\over 2} -e^{-{\hbar \over 2}} }
e^{ -{\hbar\over 4} \sum_i(f_i +{1\over f_i})} 
I_{\vec{f}, \hbar}(J_{\alpha+{1\over f}}(L)),
\end{equation}
\begin{equation}
F(U\pm)= -2{\hbar \over e^{\hbar\over 2} -e^{-{\hbar \over 2}} }
e^{-{3\over 4}\hbar \pm}
\pm.
\end{equation}
So, in conclusion~:
\begin{equation}
O(L)= e^{-{\hbar\over 4}\sum_i(f_i +{1\over f_i} -{3\over4} sgn(f_i))}
 \prod_{i=1}^{|L|}sgn(f_i)
I_{\vec{f}, \hbar}(J_{\vec{\alpha}+{1\over \vec{f}}}(L)),
\end{equation}
and if, for example, we compute the invariant for the Lens space $L(n,1)$, we 
get~:
\begin{equation}
O(L(n,1))= e^{-{ \hbar \over 4} \left(n + {2\over n} - 3 
sgn(n)\right)} { e^{\hbar \alpha \over 2|n|} - e^{-{\hbar \alpha \over 
2|n|}} \over e^{\hbar \over 2} - e^{-{\hbar \over 2}}}
\end{equation}

{\bf Acknowledgment} I would like to thank Daniel 
Altschuler for many useful discussions and his constant support
and the CNRS for its financial support.

\end{document}